\documentclass[12pt,twoside]{article}

\usepackage{graphicx,graphics, amsthm, amsfonts, amsbsy,amssymb, amsmath, cite, array,arydshln, hyperref, float,tikz}

\everymath{\displaystyle}

\usepackage[margin=1in]{geometry}

\allowdisplaybreaks

\newtheorem{theorem}{Theorem}[section]

\newtheorem{lemma}[theorem]{Lemma}
\newtheorem{corollary}[theorem]{Corollary}

\newtheorem{proposition}[theorem]{Proposition}

\newtheorem{conjecture}{Conjecture}

\begin{document}
\title{Proximity and Remoteness in Graphs: a survey}
\author{Mustapha Aouchiche$^{a}$ and Bilal Ahmad Rather$^{b,}$\footnote{Corresponding author} \\
	$^{a}$\textit{Polytechnique Montreal, Montreal, QC, Canada}\\
	$^{a}$\texttt{mustapha.aouchiche@polymtl.ca}\\
	$^{b}${\it Department of Mathematical Sciences, College of Science}, \\
	{\it United Arab Emirate University, Al Ain 15551, Abu Dhabi, UAE}\\
	$ ^{b} $\texttt{bilalahmadrr@gmail.com} 
}
\date{}

\pagestyle{myheadings} \markboth{Aouchiche and Rather}{Proximity and Remoteness in Graphs: a survey}
\maketitle
\vskip 5mm
\begin{abstract}
	The {\it proximity} $\pi = \pi (G)$ of a connected graph $G$ is the minimum, over all vertices, of the average distance from a vertex to all others. Similarly, the maximum is called the {\it remoteness} and denoted by $\rho = \rho (G)$. The concepts of proximity and remoteness, first defined in 2006, attracted the attention of several researchers in Graph Theory. Their investigation led to a considerable number of publications.  In this paper
	we present a survey of the research work done to date.
\end{abstract}
\vskip 3mm



\section{Introduction}

Models involving paths, distances and location on graphs are much studied in operations research and mathematics. Models from operations research (see e.g. \cite{laporte1,laporte2}) usually use weighted graphs to describe some well-defined class of problems or some specific applications. Models from mathematics most often consider unweighted graphs and relations between graph invariants, that is, numerical quantities whose values do not depend on the labelling of edges or vertices. In this paper, we present a survey of two  graph invariants: proximity and remoteness, defined as the minimum and maximum of the average distance from a vertex to all others.

\medskip
Let $G=(V,E)$ denote a simple and connected graph, with vertex set $V$ and edge set $E$, containing $n=|V|$ vertices and $m=|E|$ edges. All the graphs considered in the present paper are finite, simple and connected. 

For a vertex $u \in V$, the set of its \textit{neighbors} in $G$ is denoted as
$$
N_u = \{v \in V : uv \in E\}.
$$
The \textit{degree} of $u$ is the number of its neighbors, that is, $d(u) = d_u = |N_u|$. The \textit{maximum degree} $\Delta$ and the \textit{minimum degree} $\delta$ of $G$ are defined as
$$
\Delta = \max_{u \in V} d(u) \qquad \mbox{ and } \qquad \delta = \min_{u \in V} d(u).
$$
A graph $G$ is said to be $d$-\textit{regular}, or \textit{regular of degree $d$}, if all of its vertices have degree $d$.

The distance between two vertices $u$ and $v$ in $G$, denoted by $d(u,v)$, is the length of a shortest path between $u$ and $v$. The average distance between all pairs of vertices in $G$ is denoted by $\overline{l}$. The eccentricity $e(v)$ of a vertex $v$ in $G$ is the largest distance from $v$ to another vertex of $G$. The minimum eccentricity in $G$, denoted by $r$, is the radius of $G$. The maximum eccentricity of $G$, denoted by $D$, is the diameter of $G$. The \textit{girth} $g$ of a graph $G$ is the length of its smallest cycle. The average eccentricity of $G$ is denoted $ecc$. That is 
$$
r = \min_{v\in V} e(v), \qquad D = \max_{v\in V} e(v)  \qquad \mbox{ and } \qquad ecc = \frac{1}{n} \sum_{v\in V} e(v).
$$ 

\begin{figure}[h]
	\centerline{\scalebox{.35}{\includegraphics{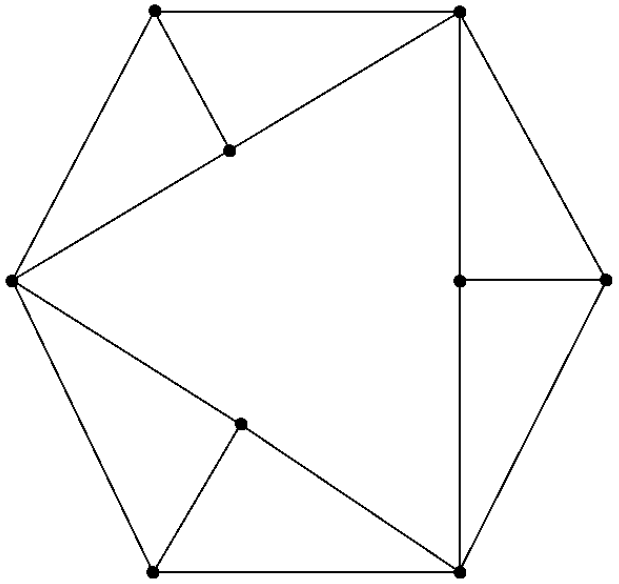}}}
	\caption{The transmission regular but not degree regular graph with the smallest order}
	\label{transreg}
\end{figure}

The {\it Wiener index} $W= W(G)$, introduced in \cite{Wiener1947}, of a connected graph $G$ is defined to be the sum of all distances in $G$, that is,
$$
W(G) = \frac{1}{2} \sum_{u, v \in V} d(u, v).
$$
The {\it transmission} $t(v)$ of a vertex $v$ is defined to be the sum of the distances from $v$ to all other vertices in $G$, that is,
$$
t(v) = \sum_{u \in V} d(u,v).
$$
A connected graph $G = (V, E)$ is said to be {\it $k$--transmission regular} if $t(v) = k$ for every vertex $v \in V$. The transmission regular graphs are exactly the {\it distance--balanced} graphs introduced in \cite{handa}. They are also called {\it self--median} graphs in \cite{cabello}. It is clear that any vertex--transitive graph (a graph $G$ in which for every two vertices $u$ and $v$, there exist an automorphism $f$ which take a vertex $u$ and map it to a vertex $v$ in $G$.
The converse is not true in general. Indeed, the graph on $9$ vertices illustrated in Figure~\ref{transreg} is  $14$--transmission regular graph but not degree regular and therefore not vertex--transitive. Actually, the graph in Figure~\ref{transreg} is the smallest transmission regular but not degree regular (see e.g. \cite{vns8,Ilic2010}). For more examples of transmission regular but not degree regular graphs see \cite{vns8,conjszeged,Ilic2010,Jerebic2008,Tavakoli2013}.

The {\it proximity } $\pi = \pi(G)$ of $G$ is the minimum average distance from a vertex of $G$ to all others. Similarly, the {\it remoteness} $\rho = \rho(G)$ of $G$ is the maximum average distance from a vertex to all others. The two last concepts were recently introduced in \cite{these,vns20}. They are close to the concept of transmission $t(v)$ of a vertex $v$. That is, if we denote $\tilde{t}(v)$ the average distance from a vertex $v$ to all other vertices in $G$, we have
$$
\pi = \min_{v\in V} \tilde{t}(v)= \min_{v\in V} \frac{t(v)}{n-1}   \qquad \mbox{ and } \qquad \rho = \max_{v\in V} \tilde{t}(v) = \max_{v \in V} \frac{t(v)}{n-1}.
$$

The transmission of a vertex is also known as the {\it distance} of a vertex \cite{entringer1976} and the minimum distance (transmission) of a vertex is studied in \cite{polansky1986}. A notion very close to the average distance from a vertex is the {\it vertex deviation} introduced by Zelinka \cite{zelinka1968} as 
$$
m_1(v) = \frac{1}{n}\sum_{u \in V} d(u,v) = \frac{t(v)}{n}.
$$

\begin{figure}
	{\centerline{\scalebox{.35}{\includegraphics{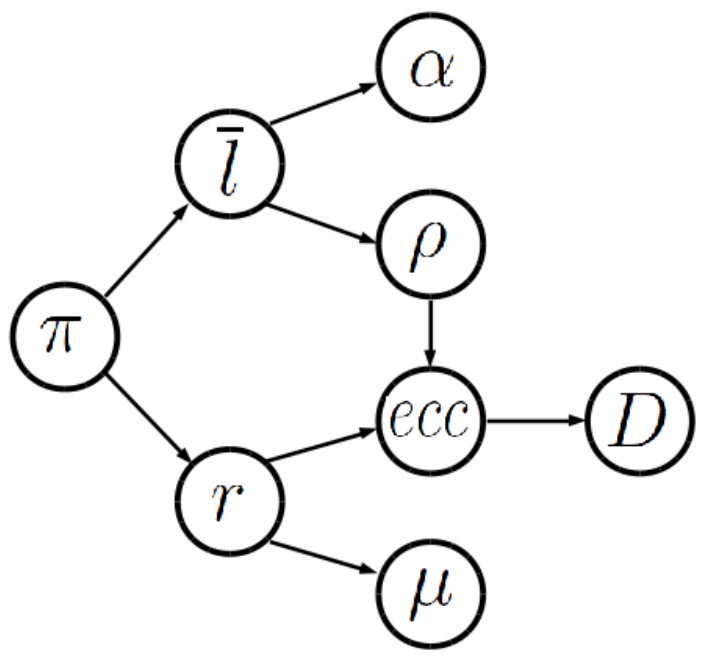}}}}
	\caption{Relations between the invariants.}
	\label{scheme}
\end{figure}

The vector composed of the vertex transmissions in a graph was first introduced by Harary \cite{harary1959} in 1959, under the name {\it the status} of a graph, as a measure of the ``weights" of individuals in social networks. The same vector was called {\it the distance degree sequence} by Bloom, Kennedy and Quintas \cite{bloom1980}. It was used to tackle the problem of graph isomorphism. Randi\'c \cite{randic1979} conjectured that two graphs are isomorphic if and only if they have the same distance degree sequence. The conjecture was refuted by several authors such as Slater \cite{slater1982}, Buckley and Harary \cite{buckley1990}, and Entringer, Jackson and Snyder \cite{entringer1976}. The transmission of a graph was also introduced by Sabidussi \cite{sabidussi} in 1966 as a measure of {\it centrality} in social networks. The notion of centrality is widely used in different branches of sciences (see for example \cite{klein2010} and the references therein) such as transportation--network theory, communication network theory, electrical circuits theory, psychology, sociology, geography, game theory and computer science. Notions closely related to that of the distance from a vertex are those of {\it a center} and {\it a centroid} introduced by Jordan \cite{jordan1869} in 1869. For mathematical properties of these two concepts see the survey, as well as the references therein, \cite{slater1999}. In 1964, Hakimi \cite{hakimi1964} used for the first time the sum of distances in solving facility location problems. In fact, Hakimi \cite{hakimi1964} considered two problems, subsequently considered in many works: the first problem was to determine a vertex $u \in V$ so as to minimize $max_{v \in V} \{ d(u,v): u \in V\}$, that is, the center of a graph; and the second problem is to determine a vertex $u \in V$ so as to minimize the sum of distances from $u$, that is, the centroid. Interpretations of these problems can be found, for instance, in \cite{goldman1972}. 
In view of the interest of the transmission vector in different domains of sciences, it is natural to study the properties of its extremal values themselves, and among the set of graph parameters. The  study of proximity and remoteness, since closely related to respectively the minimum and maximum values of the transmissions, appears to be convenient, specially, with other metric invariants, such as the diameter, radius, average eccentricity and average distance. Indeed, it follows from the definitions that 
$$
\pi \le r \le ecc \le D, \qquad \pi \le \overline{l} \le \rho \le D \qquad \mbox{ and } \qquad \overline{l} = \frac{1}{n(n-1)} \sum_{v\in V} t(v).
$$
In addition to these inequalities, several related ones can be found in the graph theory literature. Recall that a subset $S $ of vertices of $G$ is said to be {\it independent} if its vertices are pairwise non adjacent. The maximum cardinality of such a subset is called the {\it independence number} of $G$ and is denoted by $\alpha = \alpha (G)$. Then $\overline{l} \le \alpha$ \cite{chung} and $r \le \alpha$ \cite{Fajtlowicz1987}. Recall, also, that a {\it matching} in a graph is a set of disjoint edges. The maximum possible cardinality of a matching in a graph $G$ is called {\it the matching number} of $G$ and denoted by $\mu = \mu (G)$. The inequality $r \le \mu$ can be found in \cite{Fajtlowicz1987}. It was proved in \cite{Aouchiche2011} that $\rho \le ecc$. All these inequalities are illustrated in Figure~\ref{scheme}, where vertices correspond to invariants $a, b, \ldots$, and directed arcs to the relations $a \le b$. Observe that all relations mentioned are tight as all of them but $r \le \mu$ become equalities for the complete graph $K_n$ (all invariants but $\mu$ being equal to 1) and $r = \mu = 1$ for the star $S_n$. 

Since their introduction in \cite{these,vns20}, the proximity and the remoteness attracted the attention of several authors.

\section{Proximity and Remoteness}

As for any graph invariant, the first questions asked about the proximity and the remoteness are: ``what are their minimum and maximum values for given order $n$?'' and ``which extremal graphs are associated with these extremal values for a given order $n$?'' Both questions are answered in the following proposition proved in \cite{Aouchiche2011}.

\begin{proposition}[\cite{Aouchiche2011}]  
	Let $G$ be a connected graph on $n \ge 3 $ vertices with proximity $\pi$ and remoteness $\rho$. Then 
	$$
	1 \le \pi \le \left\{ 
	\begin{array}{ccl}
		\frac{n+1}{4} & & \mbox{ if $n$ is odd} \\  
		\\
		\frac{n+1}{4} + \frac{1}{4(n-1)} & & \mbox{ if $n$ is even} 
	\end{array}
	\right. 
	$$
	and 
	$$
	1 \le \rho \le \frac{n}{2}.
	$$
	The lower bound on $\pi$ is reached if and only if $G$ contains a dominating vertex; the upper bound on $\pi$ is
	attained if and only if $G$ is either the cycle $C_n$ or the path $P_n$; the lower bound on $\rho$ is reached if and
	only if $G$ is the complete graph $K_n$; the upper bound on $\rho$ is attained if and only if $G$ is the path $P_n$. 
\end{proposition}

Since the proximity $\pi$ and the remoteness $\rho$ are respectively the minimum and the maximum of the same function on a graph, it is natural to ask about how large can the difference $\rho - \pi$ be, or in other words, how large can the spread of the average distances from vertices be. This problem was solved in \cite{Aouchiche2011}.

\begin{proposition}[\cite{Aouchiche2011}]   
	Let $G$ be a connected graph on $n \ge 3$ vertices with remoteness $\rho$ and
	proximity $\pi$. Then
	$$
	\rho - \pi \le 
	\left\{
	\begin{array}{lll}
		\frac{n - 1}{4} & & \mbox{if $n $ is odd, } \\
		\frac{n-1}{4} - \frac{1}{4n - 4} & & \mbox{if $n $ is even. }
	\end{array}
	\right.
	$$
	Equality holds if and only if $G$ is a graph obtained from a path $P_{\left\lceil\frac{n}{2}\right\rceil}$ and any connected graph $H$ on $\left\lfloor n/2 \right\rfloor + 1$ vertices by a coalescence of an endpoint of the path with any vertex of $H$. 
\end{proposition}

\begin{figure}[h]
	\centerline{\scalebox{.5}{\includegraphics{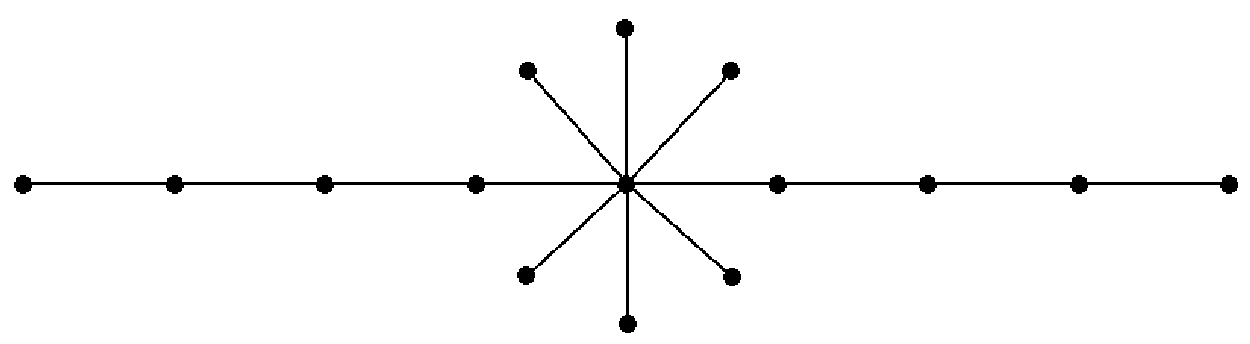}}}
	\caption{The double-tailed comet $DTC_{15,4,4}$.}
	\label{DTC}
\end{figure}

The problem of find bounds on $\pi$ and $\rho$ for a graph $G$ with given order $n$ and diameter $D$ was considered in \cite{Aouchiche2011}. Actually, the best possible lower bound on $\pi(G)$ and the best possible upper bound on $\rho(G)$ were established. These involve a particular class of graphs, next described. A {\it double-tailed comet} $DTC_{n,p,q}$ (see Figure~\ref{DTC} for an example, that is, $DTC_{15,4,4}$), with $n \ge p+q+1$, $p \ge 2$ and $q \ge 2$ is the tree obtained from a path $u_0 u_1 \cdots u_p u_{p+1} \cdots u_{p+q}$ by attaching $n-p-q-1$ pendant edges to $u_p$. It is said to be balanced if $p=q$. It is easy to see that the diameter of $DTC_{n,p,q}$ is $D = p+q$ and its maximum degree is $\Delta = d(u_p) = n-p-q+1$. Assuming, without loss of generality, that $p \ge q$, we have
\begin{eqnarray}
	\pi (DTC_{n,p,q}) & = & \frac{p(p+1)}{2(n-1)} + \frac{q(q+1)}{2(n-1)} - \frac{p+q}{n-1} +1 \nonumber \\
	& = & \frac{p(p-1)}{2(n-1)} + \frac{q(q-1)}{2(n-1)} +1 \nonumber
\end{eqnarray}
and 
$$
\rho (DTC_{n,p,q}) = \frac{(p+q)(p+q+1)}{2(n-1)} - \frac{p^2+pq}{n-1}+ p.
$$

The following result gives the lower bound for $\pi$ in terms of diameter and order of $G.$
\begin{proposition}[\cite{Aouchiche2011}]\label{diamfixe}   
	Let $G$ be a connected graph on $n$ vertices with diameter $D$. Then
	$$
	\pi (G) \ge \frac{p(p-1)}{2(n-1)} + \frac{q(q-1)}{2(n-1)} + 1
	$$
	where $p = \left\lceil D/2 \right\rceil$ and $q = \left\lfloor D/2 \right\rfloor$. The bound is best possible as shown by the double-tailed comet $DTC_{n,p,q}$.
\end{proposition}

A {\it comet} $CO_{n,\Delta}$ is the tree obtained from a star $S_{\Delta +1}$ and a path $P_{n-\Delta}$ by coalescence of an endpoint of $P_{n-\Delta}$ with a pendant vertex of $S_{\Delta +1}$. Easy computations lead to the following expressions for the diameter and the remoteness of a comet:
$$
D(CO_{n,\Delta}) = n-\Delta +1 
$$
and 
$$
\rho (CO_{n,\Delta}) = \frac{(n-\Delta +1)(n+\Delta-2)}{2(n-1)}.
$$
~\\
A {\it kite} $KI_{n,\omega}$ is the connected graph obtained from a clique $K_{\omega}$ and a path $P_{n-\omega}$ by adding an edge between an endpoint of the path and a vertex from the clique. A {\it pseudo-kite} $PKI_{n,p}$ is any connected graph which is a spanning subgraph of the kite $KI_{n,p}$ and that contains the comet $CO_{n,p}$ as a spanning tree. Note that $CO_{n,p}$, $KI_{n,p}$ and $PKI_{n,p}$ have the same proximity and the same remoteness. 

The first relationship between the pseudo-kites and the notion of remoteness is given in the following proposition.

\begin{proposition}[\cite{Aouchiche2011}]\label{degmaxfixe}   
	Let $G$ be a connected graph on $n \ge 3$ vertices with diameter $D$. Then
	$$
	\rho (G) \le \rho(PKI_{n,n-D+1})
	$$
	with equality if and only if $G$ is a pseudo-kite $PKI_{n,n-D+1}$.
\end{proposition}

Let $G$ be a graph and $\overline{G}$ its complement. If $I$ is an invariant of $G$, we denote by $\overline{I}$ the same
invariant but in $\overline{G}$. Nordhaus-Gaddum relations for the graph invariant $I$ are inequalities of the following
form
$$
L_1(n) \le I + \overline{I} \le U_1(n) \qquad \mbox{ and } \qquad L_2(n) \le I \cdot \overline{I} \le U_2(n), 
$$
where $L_1(n)$ and $L_2(n)$ are lower bounding functions of the order $n$, and $U_1(n)$ and $U_2(n)$ upper bounding
functions of the order $n$. Note that sometimes, in addition to the order $n$, other graph invariants are used in the
bounds. This type of relations is so called after Nordhaus and Gaddum \cite{norgad} who were the first authors to give
such relations, namely
\begin{equation}\label{original}
	2 \sqrt{n} \le \chi + \overline{\chi} \le n+1 \qquad \mbox{ and } \qquad n \le \chi \cdot \overline{\chi} \le
	\left(\frac{n+1}{2}\right)^2,
\end{equation}
where $\chi$ is the chromatic number of a graph. Since then many graph theorists have been interested in finding such relations for
various graph invariants. See \cite{NGSurvey} for a survey of Nordhaus--Gaddum type results. These type of inequalities relate graph invariants with their complements and information about graphs structure is encoded in some cases as it happens with some invariant, a class of graph is uniquely identified by Nordhaus--Gaddum inequality.
 For proximity and remoteness, Nordhaus--Gaddum inequalities were proved in \cite{Aouchiche2010} and are stated as: 

\begin{theorem}[\cite{Aouchiche2010}] 
	For any connected graph $G$ on $n \ge 5$ vertices for which $\overline{G}$ is connected
	$$
	\frac{2n}{n-1} \le \pi + \overline{\pi} \le \left\{ 
	\begin{array}{ccl}
		\frac{n+1}{4} + \frac{n+1}{n-1} & & \mbox{ if $n$ is odd, } \\  
		\frac{n}{4} + \frac{n}{4(n-1)} + \frac{n+1}{n-1} & & \mbox{ if $n$ is even. } 
	\end{array}
	\right. 
	$$
	The lower bound is attained if and only if $\Delta(G) = \Delta(\overline{G}) = n-2$. The upper bound is attained if and
	only if either $G$ or $\overline{G}$ is the cycle $C_n$.
\end{theorem}

\begin{theorem}[\cite{Aouchiche2010}] 
	For any connected graph $G$ on $n \ge 5$ vertices for which $\overline{G}$ is connected
	$$
	\frac{n^2}{(n-1)^2} \le \pi \cdot \overline{\pi} \le \left\{ 
	\begin{array}{ccl}
		\frac{(n+1)^2}{4(n-1)} & & \mbox{ if $n$ is odd, } \\  
		\frac{n(n+1)}{4(n-1)} + \frac{n(n+1)}{4(n-1)^2} & & \mbox{ if $n$ is even. } 
	\end{array}
	\right. 
	$$
	The lower bound is attained if and only if $\Delta(G) = \Delta(\overline{G}) = n-2$. The upper bound is attained if and
	only if either $G$ or $\overline{G}$ is the cycle $C_n$.
\end{theorem}

Recall that a {\it comet} $Co_{n,\Delta}$ is obtained from a star $S_{\Delta+1}$
by appending a path $P_{n-\Delta-1}$ to one of its pending vertices. Moreover, a {\em path-complete} graph
$PK_{n,m}$ on $n$ vertices and $m$ edges is the graph obtained from a path $P_k$, $k \ge 1$, and a clique $K_{n-k}$ by
adding at least one edge between one endpoint of the path and the vertices of $K_{n-k}$, where $ (n-k)(n-k-1)/2 + k \le
m \le (n-k+1)(n-k)/2 + k -1$. One can verify that there is exactly one path-complete graph $PK_{n,m}$ for all $n$ and
$m$ such that $1 \le n-1 \le m \le n(n-1)/2$.  

\begin{figure}[h]
	\centerline{\includegraphics[width=3.25cm]{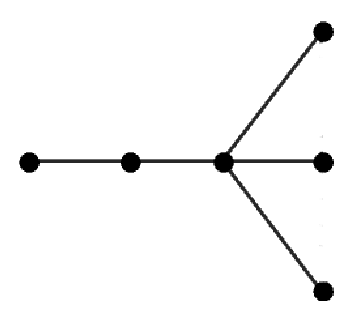}\includegraphics[width=3.25cm]{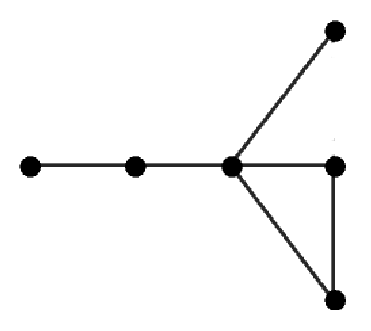}\includegraphics[width=3.25cm]{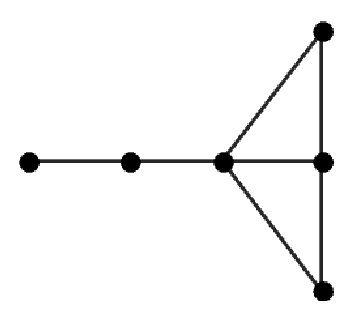}\includegraphics[width=3.25cm]{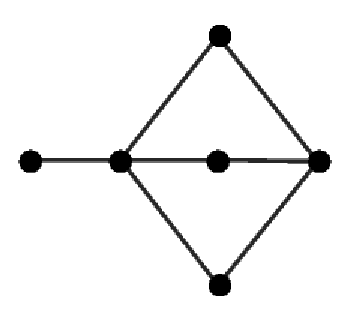}}
	\caption{Graphs with $D=3$ that maximize $\rho + \overline{\rho}$ for $n=6$.}
	\label{fig1}
\end{figure}

\begin{theorem}[\cite{Aouchiche2010}]\label{rhosum} 
	For any connected graph $G$ on $n \ge 6$ vertices for which $\overline{G}$ is connected
	$$
	3 \le \rho + \overline{\rho} \le \frac{n+2}{2} + \frac{2}{n-1}.
	$$
	The lower bound is attained if and only if $n\ge8$, $G$ is regular and $D = \overline{D} = 2$. The upper
	bound is attained if and only if $G$ or $\overline{G}$ is the path $P_n$, the comet $Co_{n,3}$ or the path-complete
	graph $PK_{n,n}$ when $n \ge 7$, and if and only if $G$ or $\overline{G}$ is the path $P_6$, the comet $Co_{6,3}$, the
	path-complete graph $PK_{6,6}$ or one of the graphs in Fig.~\ref{fig1}.  
\end{theorem}

\begin{theorem}[\cite{Aouchiche2010}]\label{thprodrho} 
	For any connected graph $G$ on $n \ge 7$ vertices for which $\overline{G}$ is connected
	$$
	\rho \cdot \overline{\rho} \le \left\{
	\begin{array}{lll}
		\frac{16n+20}{27} + \frac{8}{9(n-1)} + \frac{4}{27(n-1)^2} &  \mbox{ if } & n = \, \, 0 \,\, (mod \, \, 3), \\
		\frac{16n+20}{27} + \frac{2}{3(n-1)} &  \mbox{ if } & n = \, \, 1 \,\, (mod \, \, 3), \\
		\frac{16n+20}{27} + \frac{8}{9(n-1)} + \frac{5}{27(n-1)^2}  &  \mbox{ if } & n = \, \, 2 \,\, (mod \, \, 3). \\
	\end{array}
	\right.
	$$
	The upper bound is the best possible as shown by the comets $Co_{n, \left\lceil \frac{n}{3} \right\rceil+1}$, and 
	$Co_{n, \left\lceil \frac{n}{3} \right\rceil}$ if $n = \, \, 1 \,\, (mod \, \, 3)$.
\end{theorem}
Note that the bound provided in Theorem~\ref{thprodrho} is not valid for $n= 4, 5, 6$ (if $n \le 3$,
then $G$ and $\overline{G}$ cannot be connected simultaneously). For the path $P_4$ (respectively the comets $Co_{5,3}$
and $Co_{6,4}$), $\rho \cdot \overline{\rho} =4$ (respectively 9/2 and 24/5) while the corresponding value provided by
Theorem~\ref{thprodrho} is $10/3$ (respectively 63/16 and 112/25).

Next, we consider the some results related to the proximity for some special classes of graphs. Czabarka et al. \cite{czabarka2020} obtained results on proximity in triangulations and quadrangulations, that is, graphs with maximal planer graphs and graphs with maximal bipartite planer graphs, respectively.
\begin{theorem}[\cite{czabarka2020}]
	Let $ G $ be a planar graph of order $ n $ and $ v $ a central vertex of $ G $.
	\begin{itemize}
		\item[\bf (i)] If $ G $ is a triangulation, then
		\[ \pi \leq \frac{n + 19}{12}+\frac{25}{3(n -1)}. \]
		\item[\bf (ii)] If $ G $ is a $ 4 $-connected triangulation, then
		\[ \pi \leq \frac{n + 35}{16}+\frac{91}{4(n-1)}. \]
		\item[\bf (iii)] If $ G $ is a $ 5 $-connected triangulation, then
		\[ \pi\leq \frac{n + 57}{20}+\frac{393}{10(n -1)}. \]
		\item[\bf (iv)] If $ G $ is a quadrangulation, then
		\[ \pi \leq \frac{n + 11}{8}+\frac{9}{2(n -1)}. \]
		\item[\bf (v)] If $ G $ is a $ 3 $-connected quadrangulation, then
		\[ \pi \leq \frac{n + 25}{12}+\frac{169}{12(n-1)}. \]
	\end{itemize}
\end{theorem}

Pei et al. \cite{pei2021} gave following interesting results relating proximity, average eccentricity and domination, results earlier conjectured in \cite{these}.	
A non empty set $ S\subseteq V(G) $ is said to be a \emph{dominating} set if every vertex in $ V\setminus S $ is adjacent to at least one  vertex in $ S. $ The domination number $ \gamma (G) $ is the minimum cardinality of a dominating set of $G$.

Recall that $ \gamma(G) \leq \left \lfloor \frac{n}{2}\right \rfloor. $ In fact, the graphs with domination number $\left \lfloor \frac{n}{2}\right \rfloor$ have been determined in the following result.
\begin{figure}[h]
	\centerline{\scalebox{.98}{\includegraphics{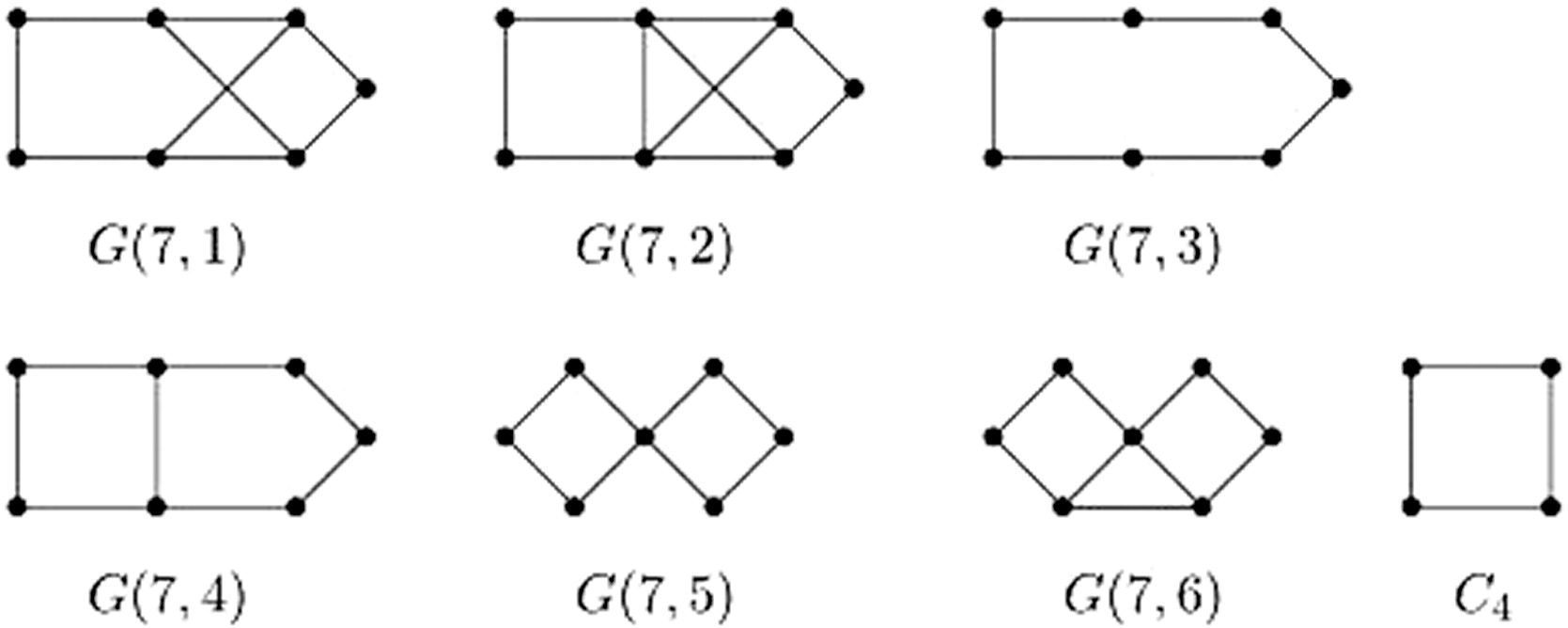}}}
	\caption{Graphs in family $\mathcal{A} $}
	\label{fancyA}
\end{figure}
\begin{figure}[h]
	\centerline{\scalebox{.98}{\includegraphics{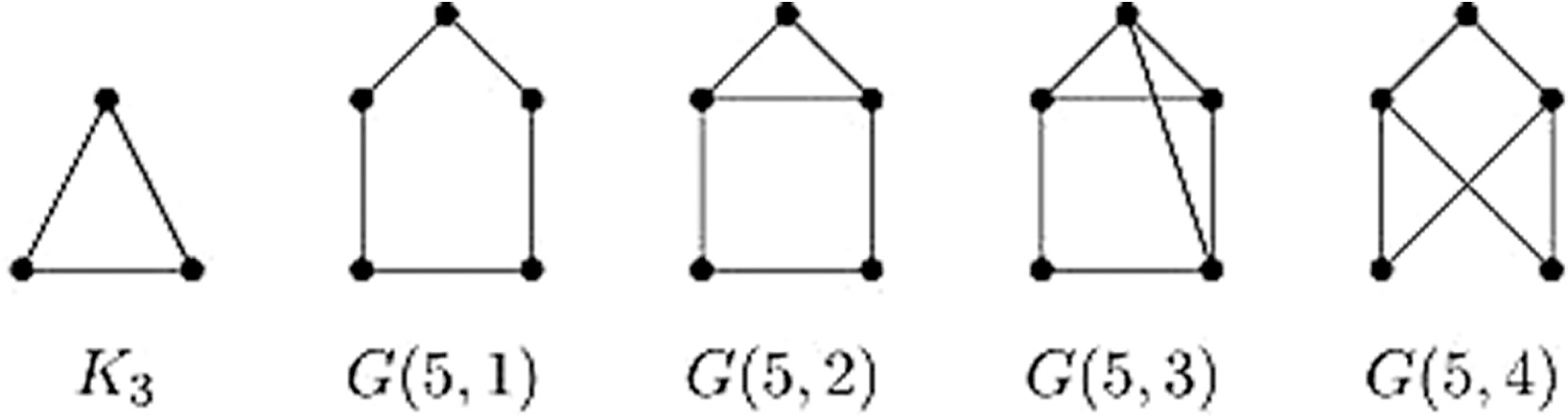}}}
	\caption{Graphs in family $\mathcal{B} $}
	\label{fancyB}
\end{figure}

\begin{lemma}[\cite{payen}]
	A connected graph $ G $ of order $ n $ satisfies $ \gamma(G) = \left \lfloor \frac{n}{2}\right \rfloor $ if and only if $ G\in \mathcal{G}= \bigcup_{i=1}^{6} \mathcal{G}_{i} $, where $ \mathcal{G}_{i},~i = 1, \dots , 6, $ is the set defined in the following.
	Let $ H $ be any graph with vertex set $ \{v_{1},\dots, v_{k}\} $. Denote by $ f(H) $ the graph obtained from $ H $ by adding new vertices $ u_{1}, \dots, u_{k} $ and the edges $ v_{i}u_{i},~ i = 1, \dots, k. $ Define
	\[ \mathcal{G}_{1 }= \{C_{4}\} \cup \{G | G = f(H)~\text{for some connected graph}~H\}. \]
	Let $\mathcal{F} = \mathcal{A} \cup \mathcal{B} $ and $ \mathcal{G}_{2} =\mathcal{ F} -\{C_{4}\}, $ where $ \mathcal{A} = \{C_{4}, G(7, i) | i = 1, \dots , 6\} $ and $ \mathcal{B} = \{K3, G(5, i) | i = 1, \dots , 4\}, $ as shown in Fig. \ref{fancyA} and Fig. \ref{fancyB}, respectively. For any graph $ H $, let $ \phi (H) $ be the set of connected graphs, each of which can be formed from $ f (H) $ by adding a new vertex $ x $ and edges joining $ x $ to one or more vertices of $ H $. Then define
	$ \mathcal{G}_{3} = \{G | G = \phi(H)~ \text{for some graph}~ H\}. $
	Let $ G \in \mathcal{G}_{3} $ and $ y $ be a vertex of a copy of $ C_{4} $. Denote by $ \theta(G) $ the graph obtained by joining $ G $ to $ C_{4} $ with the single edge $ xy, $ where $ x $ is the new vertex added in forming G. Then define
	\[ \mathcal{G}_{4} = \{G | G = \theta(H)~\text{for some graph H}~ \in \mathcal{G}_{3}\}. \]
	Let $ u, v,w $ be the vertex sequence of a path $ P_{3}. $ For any graph $ H, $ let $ \mathcal{P}(H) $ be the set of connected graphs which may be formed from $ f (H) $ by joining each of $ u $ and $ w $ to one or more vertices of M Then define
	\[ \mathcal{G}_{5} = \{G | G = \mathcal{P}(H)~\text{for some graph}~H\}. \]
	For a graph $ X \in \mathcal{B} $, let $ U \subset V(X) $ be a set of vertices such that no fewer than $ \gamma (X) $ vertices of $ X $ dominate $ V(X)\setminus U. $ Let $ \mathcal{R}(H, X) $ be the set of connected graphs which may be formed from $ f(H) $ by joining each vertex of $ U $ to one or more vertices of $ H $ for some set $ U $ as defined above and any graph $ H. $ Then define
	\[ \mathcal{G}_{6} = \{G | G \in \mathcal{R}(H, X)~\text{for some}~ X \in \mathcal{B}~ \text{and some}~ H\}. \]
\end{lemma}

\begin{lemma}[\cite{pei2021}]
	Let $ G $ be a connected graph of order $ n $ with $\gamma(G)\geq 2, \Delta(G) \leq  n-2 $ and $\delta(G) = 1. $ Then 
	\[ ecc \geq  2 + \frac{\gamma(G)}{n}. \]
\end{lemma}

\begin{lemma}[\cite{pei2021}]
	Let $ G $ be a connected graph with $ \Delta(G) \leq  n-2 $ and $\delta(G) = 2. $ Then $ \gamma(G) =2$ or $ ecc\geq 2+\frac{2}{n}. $
\end{lemma}

\begin{lemma}[\cite{pei2021}]
	Let $ G $ be a connected graph of order $ n\geq 2 $.  If $ \Delta(G) = n-1, $ then  $\frac{\gamma(G)}{ecc}\leq 1 $ with equality if and only if $ G\cong K_{n}. $
\end{lemma}

\begin{lemma}[\cite{pei2021}]
	Let $ G $ be a connected graph on $ n $ vertices, where $ 2\leq n \leq 5. $ Then 
	\[ \frac{\gamma(G)}{ecc} \leq1, \]
	with equality if and only if $ G \in \{P_{2}, C_{3}, K_{4}, C_{4}\} $ when $ n \leq 4 $, and $$ G \in  {K_{5}, G(5, 1), G(5, 2), G(5, 3), G(5, 4)} \subseteq \mathcal{B} $$ when $ n = 5, $ see Fig. \ref{fancyB}
\end{lemma}
Based on the above lemmas, the following result were presented in \cite{pei2021}.
\begin{theorem}[\cite{pei2021}]
	Let $ G $ be a connected graph of order $ n\geq 2. $ Then
	\[ 
	\frac{\gamma(G)}{ecc}\leq \begin{cases}
		1 & \text{if} ~ n\leq 5,\\
		\frac{n.\left \lfloor \frac{n}{2}\right \rfloor} {\left \lfloor \frac{5n}{2}\right \rfloor} & \text{if}~ n\geq 6,
	\end{cases} \]
	with equality if and only if $ G \in \{P_{2}, C_{3}, K_{4}, C_{4}, K_{5}, G(5, 1), G(5, 2), G(5, 3), G(5, 4)\} $ when $ n \leq  5 $, or $ G \in \{K_{\left \lceil \frac{n}{2}\right \rceil,\left \lfloor \frac{n}{2}\right \rfloor}, G(7, 1),G(7, 2), G^{\prime}, G^{\prime\prime}\} $ when $ n \geq  6, $ where  $ G^{\prime}, G^{\prime\prime} $ are defined in \cite{pei2019}.
\end{theorem}

The first observation in \cite{pei2021} is that for a vertex $ v $ in a connected graph $ G $, the proximity is
\[ \pi(v)\geq \frac{2n-2-d(v)}{n-1}, \]
with equality if and only if $ e(v)\leq 2. $

We continue with the meaning of the notations $ x, u, v,w $ in the definition of $ \mathcal{G}_{3} $ and $ \mathcal{G}_{5} $, where $ x $ is the vertex added in forming graph $ G \in \mathcal{G}_{3}, u, v,w $ is the vertex sequence of a path $ P_{3} $ mentioned in the construction of $ G_{5}. $ Let $ H^{\ast} $ be any connected graph with vertex set $ \{v_{1}, \dots, v_{k}\} $ and $ \triangle(H^{\ast}) = k -1. $ Assume that $ d_{H^{\ast}}(v_{i_{0}} ) = \triangle(H^{\ast}) $ for some $ i_{0} \in \{1,\dots, k\}, $ that is, $ \{v_{i_{0}}v_{j} | j = 1, \dots, i_{0} - 1, i_{0} + 1, \dots, k\} \subseteq E(H^{\ast}). $
Next, we define a subfamily of $ \mathcal{G} $, namely, $ \mathcal{G}^{\ast} =\bigcup_{i\in\{1,2,3,5,6\}} \mathcal{G}_{i}^{\ast}, $ where $ \mathcal{G}_{i}^{\ast} $ for $ i = 1, 2, 3, 5, 6 $ are defined as below: 
\begin{align*}
	\mathcal{G}_{1}^{\ast} &= \{G | G = f(H^{\ast})~ \text{for some connected graph}~ H^{\ast}\}, \\
	\mathcal{G}_{2}^{\ast} &= \{G(5, 2), G(5, 3), G(5, 4), G(7, 2), G(7, 5), G(7, 6)\}, \\
	\mathcal{G}_{3}^{\ast} &= \{G | G =\phi(H^{\ast})~ \text{and}~ v_{i_{0}}x \in E(G)\}, \\
	\mathcal{G}_{5}^{\ast} &= \{G | G = \mathcal{P}(H^{\ast}), v_{i_{0}}u \in E(G)~ \text{and}~ v_{i_{0}}w \in E(G)\} \\
	\text{and} &\\
	\mathcal{G}_{6}^{\ast} &= \mathcal{G}_{6}^{1} \cup \mathcal{G}_{6}^{1} \cup \mathcal{G}_{6}^{(3,2)} \cup \mathcal{G}_{6}^{(3,3)} \cup \mathcal{G}_{6}^{(3,4)}, \\
	\text{where} &\\
	\mathcal{G}_{6}^{1} &= \{G | G = \mathcal{R}(K3, H^{\ast})~ \text{with}~ |U| = 2,~\text{and}~ v_{i_{0}} s \in E(G)~ \text{for each vertex}~ s \in U\}. \\
	\mathcal{G}_{6}^{2} &= \{G | G \in \mathcal{G}_{6}~ \text{with}~ X = K_{3},~ \text{and}~ sv_{i} \in E(G)~ \text{for some vertex}~ s \in U~ \text{and}~ i = 1, \dots , k\}. 
\end{align*}
For $ i = 2, 3, 4, $
$ \mathcal{G}_{6}^{(3,i)} = \{G | G \in \mathcal{G}_{6}~ \text{with}~ X = G(5, i),~ \text{and}~ sv_{j} \in E(G)~ \text{for}~ j = 1,\dots , k, $ some vertex $ s \in U~ \text{with}~ d_{G}(5,i)(s) = 3\}.$

The following lemma presents bounds of $\pi$ for the family $\mathcal{G}^{\ast}.$ 
\begin{lemma}[\cite{pei2021}]
	Let $ G $ be a connected graph of order $ n\geq 5 $ with $ \gamma(G) =\left \lfloor \frac{n}{2}\right \rfloor $. Then
	\[ \pi \geq \begin{cases}
		\frac{3}{2} -\frac{1}{2(n-1)}, &  n~ \text{is even}\\
		\frac{3}{2}- \frac{1}{n-1}, &  n~ \text{is odd}
	\end{cases} \]
	with equality if and only if $ G \in \mathcal{G}^{\ast}. $
\end{lemma}

\begin{theorem}[\cite{pei2021}]
	Let $ G $ be a connected graph of order $ n\geq 2 $. Then
	\[ \gamma(G)-\pi \leq \begin{cases}
		\frac{n-3}{2} +\frac{1}{2(n-1)}, & \text{if $ n $ is even}\\
		\frac{n-4}{2}-\frac{1}{n-1}, &  \text{if $ n $ is odd}
	\end{cases} \]
	with equality if and only if $ G \in \mathcal{G}^{\ast}\cup \{P_{2},P_{3}, C_{3}, P_{4}, C_{4}\}. $
\end{theorem}
Recently Pei \cite{pei2022} obtained the following sequence of results for the domination number and the remoteness of a graph.
\begin{lemma}[\cite{pei2022}]
	Let $G$ be a connected graph with $n\leq 6$	vertices and  $\gamma(G)=\left \lfloor \frac{n}{2}\right\rfloor-1.$ Then for $n\geq 4, $ we have
	\[ 	\gamma(G)-\rho(G)< \begin{cases} \frac{4}{5}, & n=6\\
		0, & 4\le n\leq 5, \end{cases}\] 
		with equality if and only if $G$ is $4$-regular when $n = 6$, and $G\cong K_{n}$ when $4\leq n\leq 5.$
\end{lemma}

\begin{lemma}[\cite{pei2022}]
	Suppose that  $G$ be a connected graph of order $n\geq 7$ with  $\gamma(G)=\left \lfloor \frac{n}{2}\right\rfloor-1.$ Then
	\[ 	\gamma(G)-\rho(G)< \begin{cases} \frac{n-5}{2}+\frac{3}{2n-2}, & n ~\text{is even}\\
		\frac{n-6}{2}+\frac{2}{n-1}, & n ~\text{is odd and}~n\geq 9\\
		\frac{n-3}{4}, & n =7. \end{cases}\] 
\end{lemma}

\begin{lemma}[\cite{pei2022}]
	Suppose that  $G$ be a connected graph of order $n$ with  $1\leq \gamma(G)=\left \lfloor \frac{n}{2}\right\rfloor-2.$ Then
	\[ 	\gamma(G)-\rho(G)< \begin{cases} \frac{n-5}{2}+\frac{3}{2n-2}, & n ~\text{is even}\\
		\frac{n-6}{2}+\frac{2}{n-1}, & n ~\text{is odd and}~n\geq 9\\
		\frac{n-3}{4}, & n ~\text{is odd and}~n\leq 7. \end{cases}\] 
	
\end{lemma}

\begin{lemma}[\cite{pei2022}]
	If $G$ is a connected graph with order $n(\geq 2)$
	and   $\gamma(G)=\left \lfloor \frac{n}{2}\right\rfloor,$ then
	\[ 	\gamma(G)-\rho(G)\leq  \begin{cases}\frac{2}{3}, & n=4\\
		\frac{n-5}{2}+\frac{3}{2n-2}, & n ~\text{is even and}~n\neq 4\\
		\frac{n-3}{4}, & n ~\text{is odd and}~n\leq 7\\
		\frac{n-6}{2}+\frac{2}{n-1}, & n ~\text{is odd and}~n\geq 9,\\
	\end{cases}\] 
	with equality if and only if $G\in \{C_{4},K_{\frac{n}{2},\frac{n}{2}}~|~n~\text{is even and}~n\neq 4\}\cup (\mathcal{G}_{2}-\{G^{5}_{7})\})\cup \{K_{\lceil \frac{n}{2}\rceil,\lfloor\frac{n}{2}\rfloor},G^{\prime},G^{\prime\prime}$ $|~n ~\text{is odd and}~ n\geq 9\},$ where  $ G^{\prime}, G^{\prime\prime} $ are defined in \cite{pei2019}.
\end{lemma}

\begin{lemma}[\cite{pei2022}]
	If $G$ is a connected graph with $n\leq 7$	vertices and $\gamma(G)=\left \lfloor \frac{n}{2}\right\rfloor-1,$ then $n\geq 4$ and
	\[ 	\gamma(G)-\rho(G)\leq \begin{cases} \frac{4}{5}, & n=6\\ 0 & 4\leq n \leq 5 \end{cases}\] 
	with equality if and only if G is $4$-regular when $n=6$, and $G \cong K_{n}$ when $4\leq n \leq 5.$
\end{lemma}

Dankelmann related proximity and remoteness with minimum degree and gave several interesting results \cite{dankelmann2015} with sharp inequalities along with the characterization of graphs attaining them.

\begin{theorem}[\cite{dankelmann2015}]\label{thm 1 dankelmann2015}
	Let $ G $ be a connected graph of order $ n $ and minimum degree $ \delta $, where $ \delta\geq 2 $. Then there exists a spanning tree $ T $ of $ G $ with 
	\[ \pi(T)\leq \frac{3n}{4(\delta + 1)}+ 3 \]
	and
	\[ \rho(T)\leq \frac{3n}{2(\delta + 1)}+\frac{7}{2}. \]
\end{theorem}
The following is an immediate consequence of Theorem \ref{thm 1 dankelmann2015}
\begin{corollary}[\cite{dankelmann2015}]\
	Let $ G $ be a connected graph of order $ n $ and minimum degree $ \delta $, where $ \delta\geq 2 $. Then
	\[ \pi(G)\leq \frac{3n}{4(\delta + 1)}+ 3 \]
	and
	\[ \rho(G)\leq \frac{3n}{2(\delta + 1)}+\frac{7}{2}. \]
\end{corollary}

Dankelmann, Jonck and Mafunda \cite{dankelmann2021}, obtained bounds for $ \pi $ and $ \rho $ in triangle-free and $ C_{4} $-free graphs in terms of order and minimum degree. Before stating the results, we need the following definitions.

Given positive integers $ n, \delta $ with $ \delta\geq 3. $ Let $ A =(1, 1, \delta-1,\delta-1, 1, 1, \delta-1, \delta-1, \dots ) $ be the infinite sequence repeating the 
$ (1, 1, \delta-1, \delta-1) $-pattern indefinitely. Define the finite sequence $ X_{n,\delta} $ by
\[ X_{n,\delta} = (1, \delta,\delta- 1,a_{1},a_{2},\dots ,a_{l(A,n-4\delta)}, \delta, r_{n,\delta}), \]
where $ r_{n,\delta}=n-3\delta-\sum_{i=1}^{l(A,n-4\delta)}a_{i}. $

The \textit{sequential sum} of graphs $ G =G_{1}+G_{2}+\dots+G_{n} $ to be the sequential join such that the vertex set $ V(G) =V(G_{1})\cup V(G_{2}) \cup \dots \cup V(G_{n}) $ and the edge set $ E(G)=E(G_{1})\cup E(G_{2})\cup\dots\cup E(G_{n})\bigcup_{i=1}^{n-1}\{uv| u \in V(G_{i}),v \in V(G_{i+1})\}. $

For a finite sequence $ X=(x_{0}, x_{1},\dots, x_{d}) $ of positive integers we define the graph $ G(X) $ by
\[ G(X) = \overline{K}_{x_{0}}+ \overline{K}_{x_{1}}+\dots + \overline{K}_{x_{d}} \]

\begin{theorem}[\cite{dankelmann2021}]
	Let $ G $ be a connected, triangle free graph of order $ n $ and minimum degree $ \delta $, where $ \delta\geq 3 $ and $ n \geq6\delta. $ Then
	\[ \rho \leq \rho(G(X_{n,\delta} )). \]
\end{theorem}

\begin{corollary}[\cite{dankelmann2021}]
	If $ G $ is a connected triangle-free graph of order $ n $ and minimum degree $ \delta\geq 3 $, then
	\[\rho \leq  2\left \lceil \frac{n-3\delta}{2\delta}\right \rceil+2-\frac{\delta}{n-1}, \]
	and this bound is sharp.
\end{corollary}

\begin{theorem}[\cite{dankelmann2021}]
	If $ G $ is a connected, triangle-free graph of order $ n $ and minimum degree $ \delta\geq 3 $, then
	\[ \pi \leq  \frac{n}{2\delta}+2-\frac{5}{2\delta}-\frac{21\delta^{2} - 8\delta-3}{2\delta(n-1)}. \]
\end{theorem}

\begin{theorem}[\cite{dankelmann2021}]
	If $ G $ is a connected, $ C_{4} $-free graph of order $ n $ and minimum degree $ \delta\geq 3 $, then
	\[ \rho \leq \frac{5}{2} \left \lfloor \frac{n}{\delta^{2}-2\left\lfloor\frac{\delta}{2}\right\rfloor+1} \right \rfloor+2.\]
\end{theorem}

\begin{theorem}[\cite{dankelmann2021}]
	If $ G $ is a connected $ C_{4} $-free graph of order $ n $ and minimum degree $ \delta\geq 3 $, then
	\[\pi\leq\frac{5}{4}\left\lfloor \frac{n}{\delta^{2}-2\left \lfloor \frac{\delta}{2}\right \rfloor+1} \right\rfloor+\frac{147}{32}.\]
\end{theorem}

The next theorem shows that for many values of $ \delta $ the bound on remoteness is close to being best possible in the sense that the ratio of the coefficients of $ n $ in the bound and  approach $ 1 $ as $ \delta $ gets large.

\begin{theorem}[\cite{dankelmann2021}]
	Let $ \delta \geq 3 $ be an integer such that $ \delta=q -1 $ for some prime power $ q $. Then there exists an infinite number of $ C_{4} $-free graphs $ G $ of minimum degree at least $ \delta $ with 
	\begin{align*}
		\rho &= \frac{5}{2}\frac{n}{\delta^{2} +3\delta +2}+ O(1),\\
		\pi &= \frac{5}{4}\frac{n}{\delta^{2} +3\delta + 2}+ O(1),
	\end{align*}
	where $ n $ is the order of $ G. $
\end{theorem}

Dankelmann and Mafunda \cite{dankelmann2021a} established results about the difference between $ \pi $ and distance parameters in triangle-free and $ C_{4} $-free graphs.

\begin{theorem}[\cite{dankelmann2021a}]
	If $ G $ is a connected, triangle-free graph of order $ n\geq 7 $ and minimum degree $ \delta\geq 3 $, then
	\[ \rho-\pi\leq \frac{n + 1}{2\delta}+4. \]
\end{theorem}

\begin{theorem}[\cite{dankelmann2021a}]
	If $ G $ is a connected, $ C_{4} $-free graph of order $ n\geq 6 $ and minimum degree $ \delta\geq 3 $, then
	\[ \rho-\pi\leq \frac{5(n + 1)}{4\left (\delta^{2}-2\left\lfloor \frac{\delta}{2} \right\rfloor+1\right ) }+\frac{101}{20}. \]
\end{theorem}

\begin{theorem}[\cite{dankelmann2022}]\label{dan1}
	Let $ G $ be a connected graph of order $ n $, minimum degree $ \delta $ and maximum degree $ \triangle $. Then there exists a spanning tree $ T $ of $ G $ with
	\[\rho(T)\leq \frac{3(n^{2}-\triangle^{2})}{2(n-1)(\delta+1)}+ 7. \]
\end{theorem}
Since $ \rho(G)\leq \rho(T) $ for every spanning tree of a connected graph $ G, $ the following result follows from Theorem \ref{dan1}.
\begin{corollary}[\cite{dankelmann2022}]
	If $ G $ is a connected graph of order $ n $, minimum degree $ \delta $ and maximum degree $ \triangle $, then
	\[\rho\leq \frac{3(n^{2}-\triangle^{2})}{2(n-1)(\delta+1)}+ 7. \]
\end{corollary}

\begin{theorem}[\cite{dankelmann2022}]
	Let $ G $ be a connected graph of order $ n, $ minimum degree $ \delta $ and maximum degree $ \triangle $.\\
	If $ \triangle> \frac{n}{2}-1, $ then 
	\[ \pi \leq \frac{3(n-\triangle)^{2}}{2(n-1)(\delta+1)}+\frac{13}{2}. \]
	If $ \triangle\leq \frac{n}{2}-1, $ then
	\[ \pi \leq \frac{3(n^{2}-2\triangle^{2})}{4(n-1)(\delta+1)}+\frac{35}{4}. \]
\end{theorem}

A thesis about the proximity and the remoteness of graphs written by Mallu in 2022 can be seen \cite{arif}, jointly supervised by Dankelmann and Mafunda.

\begin{figure}[h]
	\centerline{\scalebox{.4}{\includegraphics{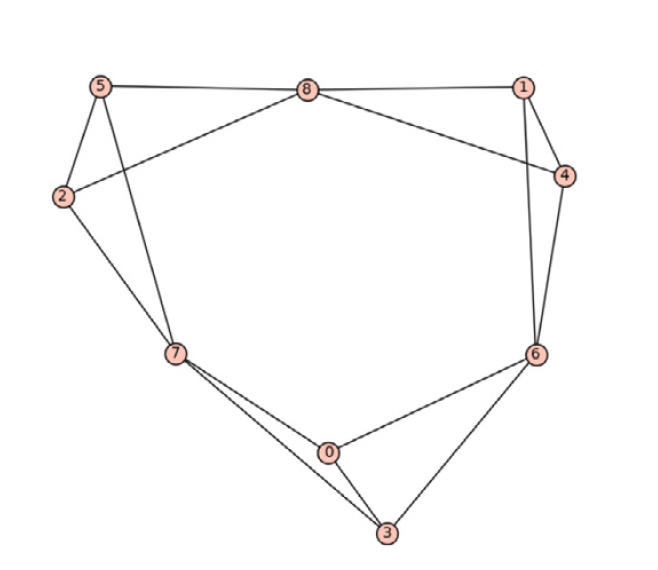}}\scalebox{.4}{\includegraphics{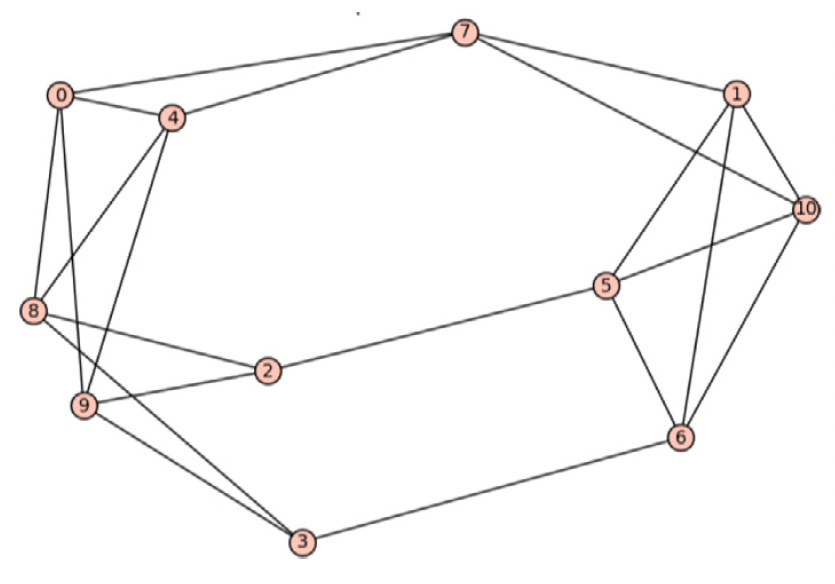}}}
	\caption{A non-regular graph $G$ of order $9$ and $11$ with $\rho(G) = \pi(G)$.}
	\label{ai1}
\end{figure}

How about $\rho(G) = \pi(G)$, where $G$ is a simple graphs? By Computer search the authors in \cite{ai} found two non-regular graphs with $\rho(G) = \pi(G)$. Figure \ref{ai1} shows two non-regular graphs of order $9$ (left) and order $11$ (right) with $\rho(G) = \pi(G)$.

The next result gives the infinite family of non-regular graphs with $\rho(G) = \pi(G).$
\begin{theorem}[\cite{ai}]
	There is an infinite number of non-regular graphs $G$ with $\rho(G) = \pi(G).$
\end{theorem}

Although, survey is primarily related the simple graphs, but for the sake of completeness and reviewer suggestions, we mention some results for digraphs.
Let $D = (V(D), A(D))$ be a strongly connected (often shortened to ‘‘strong’’) digraph with $n \geq 2$ vertices ($D$ is strong if for every pair $u, v$ of vertices in $D$ there are a path from $u$ to $v$ and a path from $v$ to $u$). We call $n = |V(D)|$ the order of $D$ and $m = |A(D)|$ the size of $D.$ The distance $d(u, v)$ from vertex $u$ to vertex $v$ in $D$ is the length of a shortest dipath from $u$ to $v.$ The distance of a vertex $u \in V(D)$ is defined as $\sigma(u) =\sum_{v\in V(D)} d(u, v),$ and the average distance of a vertex as $\overline{\sigma}(u) = \frac{\sigma(u)}{n-1}.$ The proximity $\pi(D)$ and the remoteness $\rho(D)$ are $\min{\overline{\sigma}(u)| u \in V(D)}$ and $\max{ \overline{\sigma}(u)| u \in V(D)}$, respectively. If $v$ is a vertex in a digraph $D$, then the out-neighborhood of $v$ is $N^{+}(v) = \{u \in V(D) | vu \in A(D)\}$ and in-neighborhood of $v$ is $N^{-}(v) = \{u \in V(D) | uv \in A(D)\},$ and the out-degree $d^{+}_{D}(v)$ and in-degree $d^{-}_{D}(v)$ of $v$ are $|N^{+}(v)|$ and $|N^{-}(v)|,$ respectively. An orientation of a graph $G$ is a digraph obtained from $G$ by replacing each edge by exactly one of the two possible arcs. We call an orientation of a complete graph and complete k-partite graph a tournament and a k-partite tournament,respectively. A 2-partite tournament is also called a bipartite tournament.

The very first result establishes the sharp bounds for $\pi$ and $\rho$ of $D$ along with the characterization of graphs attaining them.
\begin{theorem}[\cite{ai}]\label{ai thm 1}
	Let $D$ be a strong digraph on $n \geq 3$ vertices. Then
	\[ 1 \leq \pi(D) \leq \frac{n}{2}. \]
	The lower bound holds with equality if and only if $D$ has a vertex of out-degree $n-1.$ The upper bound holds with equality if and only if $D$ is a dicycle. And,
	\[ 1 \leq \rho(D) \leq \frac{n}{2}. \]
	The lower bound holds with equality if and only if $D$ is a complete digraph. The upper bound holds with equality if and only if $D$ is strong and contains a Hamiltonian dipath $v_{1}v_{2}\dots v_{n}$ such that $\{v_{i}v_{j} | 2 \leq i + 1 < j \leq n\} \subseteq A(D).$
\end{theorem}

From the above theorem, we have the following result.
\begin{theorem}[\cite{ai}]
	Let D be a strong digraph. Then 
	\[ \rho(D) - \pi(D) \leq \dfrac{n}{2}- 1. \]
	The upper bound holds with equality if and only if $D$ contains a Hamiltonian dipath $v_{1}v_{2}\dots v_{n}$ such that $\{v_{i}v_{j} | 2 \leq i+1 < j \leq n\} \subseteq A(D)$ and at least one of the vertices $v_{n-1}, v_{n}$ has out-degree $n-1.$
\end{theorem}

A natural question arises for digraphs: When is $\rho(D)=\pi(D)$? The answer is positive for strong tournaments. Now, we have bounds for $\pi$ and $\rho$ of tournaments, which is an analog of Theorem \ref{ai thm 1} and it show that for a strong tournament $T$, $\rho(T ) = \pi(T )$ if and only if $T$ is regular.

\begin{theorem}[\cite{ai}]\label{them 2}
	Let $D$ be a strong tournament on $n$ vertices. Then
	\[ \frac{n}{n -1}\leq  \pi(D) \leq \begin{cases}
		\qquad\frac{3}{2} & \text{if $n$ is odd}\\
		\frac{3}{2}-\frac{1}{2(n-1)} & \text{if $n$ is even.}
	\end{cases}\]
	The lower bound holds with equality if and only if $\Delta^{+}(D) = n -2.$ The upper bound holds with equality if and only if $D$ is a regular or almost regular tournament. And,
	\[ \frac{n}{2}\geq  \rho(D) \geq \begin{cases}
		\qquad\frac{3}{2} & \text{if $n$ is odd}\\
		\frac{3}{2}+\frac{1}{2(n-1)} & \text{if $n$ is even.}
	\end{cases}\]
	The lower bound holds with equality if and only if $D$ is a regular or almost regular tournament. The upper bound holds with equality if and only if $D$ is isomorphic to the tournament $T_{n}$ with $V(T_{n}) = \{v_{1},\dots,v_{n}\}$ and $A(T_{n}) = \{v_{i}v_{i+1} | i \in [n-1]\}\cup \{v_{j}v_{i} |2 \leq i + 1 < j \leq n\}.$
\end{theorem}

Theorem \ref{them 2}  allows us to easily obtain  a characterization of strong tournaments $T$ with $\rho(T ) = \pi(T )$ in the following result.
\begin{theorem}[\cite{ai}]
	For any strong tournament $T$, we have $\rho(T ) = \pi(T )$ if and only if $T$ is a strong regular tournament.
\end{theorem}

Based on the above results, one may conjecture that every digraph D with $\rho(D) = \pi(D)$ is regular. However, this is not true as such counterexamples can be found already among bipartite tournaments (see \cite{ai}). We need the following definition:

A vertex $x$ of a digraph $D$ is called bad if the out-neighborhood of $x$ is a proper subset of the out-neighborhood of another vertex of D. A vertex $z$ is good if it is not bad. A bipartite tournament is called bad if it has a bad vertex, otherwise it is called good. For a vertex $v$, we denote by $M(v)$ the set of vertices with the same out-neighborhood as $v$, that is, $M(v)=\{u|N^{+}(u) = N^{+}(v)\}$ and we denote $|M(v)|$ by $\mu(u).$ Now, we have the following results from \cite{ai}.

\begin{lemma}[\cite{ai}]
if	$T=T[A, B]$ be a strong bipartite tournament with partite sets $A$ and $B$ of sizes $n$ and $m,$ respectively, then $\pi(T )\neq \rho(T ).$
\end{lemma}

\begin{lemma}[\cite{ai}]
	Let $T = T [A, B]$ be good, $\rho(T ) =\pi(T )$ and $u, v \in V(T )$ such that $u\neq v.$ If both $u$ and $v$ are in $A$ or in $B,$ then
	\[ \mu(u) - d^{+}(u) = \mu(v) - d^{+}(v). \]
	If $v \in A$ and $u \in B$ then
	\[ 2(\mu(v) - d^{+}(v)) + |B| = 2(\mu(u) - d^{+}(u)) + |A|. \]
\end{lemma}

The above two lemmas imply the following result.
\begin{corollary}[\cite{ai}]
	For a strong bipartite tournament $T$ , we have $\rho(T ) =\pi(T )$ if and only if $T$ is good and there is a constant $c$ such that for every $v\in A$ and $u\in B,$
	\[ 2(\mu(v) - d^{+}(v)) + m = 2(\mu(u) - d^{+}(u)) + n = c. \]
	In particular, for a strong bipartite tournament $T$ with $n = m,$ we have $\rho(T ) =\pi(T )$ if and only if $T$ is good and $d^{+}(u)-\mu(u)$ is the same for every vertex $u.$
\end{corollary}

\begin{corollary}[\cite{ai}]
	Let $T$ be a good bipartite tournament and $c$ a constant such that $\mu(x) = c$ for every $x \in V(T ).$ Then $\rho(T ) = \pi(T )$ if and only if $d^{+}(v) = \frac{m}{2}$ and $d^{+}(u) = \frac{n}{2}$ for every $v\in A$ and $u\in B.$
\end{corollary}

The following result shows, in particular, that there are non-regular digraphs $D$ for which $\rho(D) = \pi(D).$
\begin{theorem}[\cite{ai}]
	For both $|A| = |B|$ and $|A|\neq |B|,$ there is an infinite number of bipartite tournaments $T$ with $\rho(T ) = \pi(T ).$
\end{theorem}

\section{Proximity and Remoteness Compared to other Metric Invariants}
Aouchiche and Hansen \cite{Aouchiche2011} in 2011 compared proximity and remoteness  with the metric invariants of a graph, most notably like diameter, radius average eccentricity and other invariants.
\begin{proposition}[\cite{Aouchiche2011}]   
	Let $G$ be a connected graph on $n \ge 3$ vertices with diameter $D $ and
	proximity $\pi $. Then
	$$
	D - \pi \le 
	\left\{
	\begin{array}{ll}
		\frac{3n-5}{4} & \quad \mbox{ if $n $ is odd, } \\
		\frac{3n-5}{4} - \frac{1}{4 n -4} & \quad \mbox{ if $n $ is even, }
	\end{array}
	\right.
	$$
	with equality if and only if $G $ is a path $P_n $.
\end{proposition}

\begin{proposition}[\cite{Aouchiche2011}]  
	Let $G $ be a connected graph on $n\ge 3$ vertices with remoteness $\rho $ and diameter $D $. Then
	$$
	D - \rho \le \frac{n-2}{2}, 
	$$
	with equality if and only if $G $ is the path $P_n$.
\end{proposition}

\begin{proposition}[\cite{Aouchiche2011}]   
	Let $G$ be a connected graph on $n \ge 3$ vertices with radius $r $ and
	proximity $\pi $. Then
	$$
	r - \pi \le 
	\left\{
	\begin{array}{ll}
		\frac{n-1}{4}-\frac{1}{4(n - 1)}  & \mbox{if $n $ is even, } \\
		\frac{n - 1}{4} - \frac{1}{n - 1} & \mbox{ if $n$ is odd. } 
	\end{array}
	\right.
	$$
	The bound is best possible as shown by the graph composed of a cycle with an additional edge forming a triangle or two additional crossed edges on four successive vertices of the cycle if $n $ is odd, and by the path $P_n$ or the cycle $C_n$ if $n$ is even.
\end{proposition}

\begin{proposition}[\cite{Aouchiche2011}]   
	Let $G$ be a connected graph on $n\ge 3$ vertices with remoteness $\rho $ and radius $r $. Then
	$$
	\rho - r \le \left\{
	\begin{array}{lll}
		\frac{n+1}{8} + \frac{1}{8(n-1)} & & \mbox{ if } n = 0 \, \, (mod \, \,  4), \\
		\frac{n+1}{8}  & & \mbox{ if } n = 1 \, \, (mod \, \,  4), \\
		\frac{n+1}{8} - \frac{3}{8(n-1)} & & \mbox{ if } n = 2 \, \, (mod \, \,  4), \\
		\frac{n+1}{8}  & & \mbox{ if } n = 3 \, \, (mod \, \,  4). \\
	\end{array}
	\right.
	$$
	The bound is attained if and only if $G$ is a pseudo-kite $PKI_{n,n-2r^*+1}$ where 
	$$
	r^* = \left\{
	\begin{array}{lll}
		\frac{n}{4} & & \mbox{ if } n = 0 \, \, (mod \, \,  4) , \\
		\frac{n-1}{4}  & & \mbox{ if } n  = 1 \, \, (mod \, \,  4), \\
		\frac{n-2}{4} \mbox{ or } \frac{n+2}{4} & & \mbox{ if } n = 2 \, \, (mod \, \,  4), \\
		\frac{n+1}{4}  & & \mbox{ if } n = 3 \, \, (mod \, \,  4). \\
	\end{array}
	\right.
	$$
\end{proposition}

\begin{proposition}[\cite{Aouchiche2011}]   
	Let $G$ be a connected graph on $n\ge 3$ vertices with remoteness $\rho $ and average eccentricity
	$ecc $. Then
	$$
	\rho \le ecc 
	$$
	with equality in both cases if and only if $G $ is a complete graph $K_n $.
\end{proposition}

Ma, Wu and Zhang \cite{ma2012} gave the following lemmas and  proved a conjecture that consists of an upper bound on $ecc - \pi$ together with the corresponding extremal graphs. 

A vertex $ v\in V $ is called a \textit{centroidal vertex} if $ \pi(v)=\pi(G) $, and the set of all centroidal vertices is the centroid (sometimes known as median or barycenter) of $ G. $

For an edge $ e \in E(G) $ and a vertex $ u \in V(G) $, we denote by $ (G-e)_{u} $ the component of $ G-e $ containing $ u $, and let $ V_{u}(e) = V((G-e)_{u}) $ and $ n_{u}(e) =|V_{u}(e)| $. Clearly, for any edge $ e = uv \in E(G),~ n_{u}(e)+n_{v} (e) = |V (G)| $. The following lemmas are given in \cite{ma2012}.
\begin{lemma}[\cite{ma2012}]\label{lem 1 ma2012}
	Let $ G $ be a tree of order $ n $. For any edge $ e = uv \in E(G), $
	\[   \pi(u) + \frac{1}{n-1} n_{u}(e) = \pi(v)+ \frac{1}{n-1}n_{v}(e). \]
\end{lemma}

\begin{lemma}[\cite{ma2012}]\label{lem 2 ma2012}
	The following holds for a tree $ G $ of order $ n $
	\begin{itemize}
		\item[\bf (i)] If $ x $ is a centroidal vertex of $ G $, then
		\[ ecc(x)\leq \left \lfloor \frac{n}{2} \right \rfloor, \]
		with equality if and only if there exists a path of length $ \left \lfloor \frac{n}{2} \right \rfloor $ in $ G $,which joins $ x $ to a pendent vertex of $ G $ with the property that the degree of every internal vertex of it is equal to two in $ G $.
		\item[\bf (ii) ] If there is a path $ P $ of length $ \left \lfloor \frac{n}{2} \right \rfloor $ in $ G $, which joins a vertex $ y $ and a pendent vertex of $ G $ with the property that the degree of every internal vertex of it is equal to two in $ G $, then $ y $ is a centroidal vertex of $ G $.
	\end{itemize}
\end{lemma}
\begin{lemma}[\cite{ma2012}]\label{lem 3 ma2012}
	Let G be a tree of order $ n \geq 3 $. Let $ v_{0} v_{1} \dots v_{d} $ be a longest path in $ G $. Set $ V_{0} = \{v \in V (G) | d(v, v_{0})\geq d(v, v_{d})\},~V_{d} = \{v \in V (G) | d(v, v_{d})\geq d(v, v_{0})\} $. Without loss of generality, let $ |V_{0}| \geq|V_{d}| $. If $ G $ is not a path, then the following holds:
	\begin{itemize}
		\item[\bf (1)] there is a pendent vertex $ v^{\ast} \in V_{0} $ distinct from $ v_{d} $,
		\item[\bf (2)] $ ecc(G)-ecc(G^{\prime})\geq \frac{1}{2} $, where $ G^{\prime} $ is the tree obtained from $ G $ by deleting the edge incident with $ v^{\ast} $ and joining $ v^{\ast} $ and $ v_{0} $,
		\item [\bf (3)] $ \pi(G)-\pi(G^{\prime}) < \frac{1}{2} $, where $ G^{\prime} $ is the tree as defined in (2).
	\end{itemize}
\end{lemma}

Lemmas \ref{lem 1 ma2012}, \ref{lem 2 ma2012} and \ref{lem 3 ma2012} help us in establishing the following result and verifies a conjecture of \cite{Aouchiche2011}.

\begin{theorem}[\cite{ma2012}]
	For a connected graph $G$ on $n \ge 3$ vertices,
	$$
	ecc - \pi \le 
	\left\{
	\begin{array}{lll}
		\frac{(3n+1)(n-1)}{4n} - \frac{n+1}{4} & & \mbox{ if $n$ is odd,}\\
		\frac{n-1}{2} - \frac{n}{4(n-1)} & & \mbox{ if $n$ is even,}
	\end{array}
	\right.
	$$
	with equality if and only if $G$ is the path $P_n$.
\end{theorem}

Sedlar \cite{sedlar2013} studied three AutoGraphiX conjectures involving proximity and remoteness. She solved the conjecture by using the following graph transformation for trees.
Let $ G $ be a tree. Then she used graph transformation (see proof in \cite{sedlar2013} of lemmas stated below) which transforms tree $ G $ to $ G^{\prime} $, where $ G^{\prime} $ is either:
\begin{itemize}
	\item[1)] Path $ P_{n} $,
	\item[2)] a tree consisting of four paths if equal length with a common end point,
	\item[3)] a tree consisting of three paths of almost equal length with a common end point.
\end{itemize}
Next, the following results prove that among those graphs the difference $ \overline{ l}-\pi $ is maximum for the last.

\begin{lemma}[\cite{sedlar2013}] \label{lem 1 sedlar2013}
	The difference $\overline{ l}-\pi$ is greater for a tree $ G $ on $ n \geq4  $ vertices consisting of three paths of almost equal length with a common end point than for path $ P_{n} $.
\end{lemma}

\begin{lemma}[\cite{sedlar2013}]\label{lem 2 sedlar2013}
	The difference $\overline{ l}-\pi$ is greater for the tree $ G $ on $ n \geq 9 $ vertices, where $ n =1 $mod$(4) $, consisting of three paths of almost equal length with a common end point than for the tree $ G^{\prime}$ on $ n $ vertices consisting of four paths of equal length, where $ G^{\prime} $ is tree obtained from $ G $ by some transformation.
\end{lemma}

\begin{lemma}[\cite{sedlar2013}] \label{lem 3 sedlar2013}
	Let $ G $ be a tree on $ n \geq 6 $ vertices with at least four leafs (a vertex of degree $ 1 $ in a tree).Then there is a tree $ G^{\prime} $ on $ n $ vertices with three leafs for which the difference $\overline{ l}-\pi$ is greater or equal than for $ G$.
\end{lemma}

\begin{lemma}[\cite{sedlar2013}] \label{lem 4 sedlar2013}
	Among trees with three leafs, the difference $\overline{ l}-\pi$ is maximal for the tree $ G $ on $ n $ vertices consisting of three paths of almost equal length with a common end vertex.
\end{lemma}

The main conclusion of Lemmas \ref{lem 1 sedlar2013}, \ref{lem 2 sedlar2013}, \ref{lem 3 sedlar2013} and \ref{lem 4 sedlar2013} can be summarized into the following result.
\begin{theorem}[\cite{sedlar2013}] \label{thm 1 sedlar2013}
	Among all trees on $ n \geq 4~ (n ,\neq 5) $ vertices with average distance $ \overline{ l} $ and proximity $ \pi $, the difference $\overline{ l}-\pi$ is maximal for a tree $ G $ composed of three paths of almost equal lengths with a common end vertex.
\end{theorem}

From Theorem \ref{thm 1 sedlar2013}, the following result follows and settles a conjecture involving proximity and remoteness.
\begin{theorem}[\cite{sedlar2013}]
	Among all connected graphs $G$ on $n \ge 3$ vertices with average distance $\overline{ l}$ and proximity $\pi$, the difference $\overline{ l} - \pi$ is maximum for a graph $G$ composed of three paths of almost equal lengths with a common end point.
\end{theorem}
Sedlar \cite{sedlar2013} also proved partial results related to another conjecture: while the conjecture is stated for all connected graphs, the results are proved for trees only.
\begin{lemma}[\cite{sedlar2013}] \label{lem 5 sedlar2013}
	Let $ G $ be a tree on $ n $ vertices with diameter $ D $ and let $ P = v_{0}v_{1}\dots v_{D} $ be a diametric path in $ G $. If there is $ j \leq \frac{D}{2} $ such that the degree of $ v_{k} $ is at most $ 2 $ for $ k \geq j+1 $, then the difference $ ecc-\rho$ is greater or equal for the path $ P_{n} $ than for $ G. $
\end{lemma}

\begin{theorem}[\cite{sedlar2013}]
	Among all trees on $n \ge 3$ vertices, the difference $ecc - \rho$ is maximum for the path $P_n$.
\end{theorem}

The following sequence of results are given by Sedlar \cite{sedlar2013}, which relates remoteness $ \rho $ with radius $ r. $
\begin{lemma}[\cite{sedlar2013}]\label{lem 5.1 sedlar2013}
	Let $ G $ be a tree on $ n $ vertices. There is a caterpillar tree $ G^{\prime} $ on $ n $ vertices for which the difference $ \rho -r $ is less or equal than for $ G $.
\end{lemma}

\begin{lemma}[\cite{sedlar2013}]\label{lem 5.2 sedlar2013}
	Let $ G \neq P_{n} $ be a caterpillar tree on $ n $ vertices with diameter $ D $,remoteness $ \rho $ and only one centroidal vertex. Let $ P = v_{0}v_{1}\dots v_{D} $ be the diametric path in $ G $ such that $ v_{j} \in P $ is the only centroidal vertex in $ G $ and every of the vertices $ v_{j+1},\dots, v_{D} $ is of the degree at most $ 2 $. Then there is a caterpillar tree $ G^{\prime}$ on $ n $ vertices of the diameter $ D+1 $ and the remoteness at most $ \rho+\frac{1}{2}. $
\end{lemma}

\begin{lemma}[\cite{sedlar2013}]\label{lem 5.3 sedlar2013}
	Let $ G \neq P_{n} $ be a caterpillar tree on $ n $ vertices with diameter $ D $, remoteness and exactly two centroidal vertices. Let $ P = v_{0}v_{1} \dots v_{D} $ be a diametric path in $ G $ such that $ v_{j},v_{j+1} \in P $ are centroidal vertices and every of the vertices $ v_{j+1},\dots,v_{D}$ is of degree at most $ 2 $. Then there is a caterpillar tree $ G^{\prime} $ on $ n $ vertices of the diameter $ D+1 $ and the remoteness at most $\rho+ \frac{1}{2}. $
\end{lemma}

\begin{lemma}[\cite{sedlar2013}]\label{lem 5.4 sedlar2013}
	Let $ G \neq P_{n} $ be a caterpillar tree on $ n $ vertices with diameter $ D $, remoteness $ \rho $ and exactly two centroidal vertices of different degrees. Let $ P=v_{0}v_{1} \dots v_{D} $ be a diametric path in $ G $ such that $ v_{j}, v_{j+1} \in P $ are centroidal vertices and every of the vertices $ v_{0},\dots, v_{j-1}, v_{j+2},\dots, v_{D} $ is of degree at most $ 2 $. Then there is a caterpillar tree $ G^{\prime} $ on $ n $ vertices of the diameter $ D +1 $ and the remoteness at most $\rho+ \frac{1}{2}  $.
\end{lemma}

\begin{lemma}[\cite{sedlar2013}]\label{lem 5.5 sedlar2013}
	Let $ G \neq P_{n}$ be a caterpillar tree on $ n $ vertices with diameter $ D $, remoteness $ \rho $ and exactly two centroidal vertices of equal degrees. Let $ P = v_{0}v_{1}\dots v_{D} $ be a diametric path in $ G $ such that $ v_{j}, v_{j+1} \in P $ are centroidal vertices and every of the vertices $ v_{0}, \dots, v_{j-1}, v_{j+2}, \dots, v_{D} $ is of degree at most $ 2 $. Then the difference $\rho-r $ is less or equal for path $ P_{n} $ than for $ G. $
\end{lemma}

\begin{lemma}[\cite{sedlar2013}]\label{lem 5.6 sedlar2013}
	Let $ G $ be a caterpillar tree on $ n $ vertices. If $ n $ is odd, then the difference $\rho-r $ is less or equal for path $ P_{n} $ then for $ G $. If $ n $ is even, then the difference $\rho- r $ is less or equal for path $ P_{n}-1 $ with a leaf appended to a central vertex than for G.
\end{lemma}

Lemmas \ref{lem 5.1 sedlar2013}, \ref{lem 5.2 sedlar2013}, \ref{lem 5.3 sedlar2013}, \ref{lem 5.4 sedlar2013}, \ref{lem 5.5 sedlar2013} and \ref{lem 5.6 sedlar2013}
can be summarized in the following result, which gives minimal trees for $ \rho-r. $
\begin{theorem}[\cite{sedlar2013}]
	Let $ G $ be a tree on $ n $ vertices. If $ n $ is odd, then the difference $\rho-r $ is less or equal for path $ P_{n} $ then for $ G $. If $ n $ is even, then the difference $\rho-r $ is less or equal for path $ P_{n-1} $ with a leaf appended to a central vertex than for $ G $.
\end{theorem}

Let $ G $ be a connected graph. A vertex $ u \in V (G) $ is called a peripheral vertex if $ \sigma(u)=\rho(G) $. For a vertex $ u\in V (G), $ let $ V_{i}(u) = \{v\in  V (G)| d(u,v) = i\} $ and $ n_{i}(u) = |V_{i}(u)| $ for each $ i \in \{1,2,\dots,d\}, $where $ d = e_{G}(u). $ In what follows, $ V_{i}(u) $ is simply denoted by $ V_{i} $ for a peripheral vertex $ u $ of $ G. $	
Wu and Zhang \cite{wu2014} proved some lemmas and two theorems, first conjectured in \cite{Aouchiche2011}.

\begin{lemma}[\cite{wu2014}]\label{lemma 1 wu2014}
	Let $ G $ be a connected graph of order $ n\geq 3. $ Let $ u $ be a peripheral vertex of $ G $ and let $ d = e_{G}(u). $ Let $ G^{\prime}$ be the graph obtained from $ G $ by joining each pair of all non adjacent vertices $ x,y $ of $ G, $ where $ x,y \in V_{j}\cup V_{j+1} $ for some $ j \in \{1,2,\dots,d-1\} $. We have
	\[ \rho(G^{\prime})-\overline{ l}(G^{\prime}) \geq \rho(G)-\overline{ l}(G),\]
	with equality if and only if $ G^{\prime} = G. $
\end{lemma}

\begin{lemma}[\cite{wu2014}]\label{lemma 2 wu2014}
	Let $ G $ be a connected graph of order $ n\geq 3 $. Let $ u $ be a peripheral vertex of $ G $ and $ e_{G}(v) = d. $ Assume that $ G[V_{j}\cup V_{j+1}] $ is a clique for each $ j\in \{0,1,\dots,d-1\} $. Let $ G^{\prime} $ be the graph with $ V (G^{\prime}) = V (G) $ and $ E(G^{\prime}) = E(G)\cup [ xy: x\in V_{d}, y\in V_{d}. $ If $ d > \left \lfloor \frac{n+1}{2}\right \rfloor, $ then
	\[ \rho(G^{\prime})-\overline{ l}(G^{\prime}) \leq \rho(G)-\overline{ l}(G), \]
	with equality if and only if $ n $ is even and $ d=\frac{n}{2}+11. $
\end{lemma}

\begin{lemma}[\cite{wu2014}]\label{lemma 3 wu2014}
	Let $ G $ be a connected graph of order $ n\geq 3 $. Let $ u $ be a peripheral vertex of $ G $ and $ e_{G}(v)=d. $ Assume that $G[V_{j}\cup V_{j+1}] $ is a clique for each $j\in \{0,1,\dots,d-1\} $. Let $ i $ be the smallest integer in $ \{1,2,\dots,d\} $ such that $ n_{i}(u)\geq 2. $ Let $ V_{i-1}(u)=\{u_{i-1}\} $ and $ v $ a vertex in $ V_{i}(u) $. Denote by $ G^{\prime} $ the graph with $ V (G^{\prime})=V (G) $ and $ E(G^{\prime})=E(G)\setminus (\{ u_{i-1}y: y\in V_{i}\setminus \{v\}\}\cup A), $ where $ A=\{vx:x\in V_{i+1}\}, $ if $ i\leq d-1 $, and $ A=\emptyset$ otherwise. If $ d<\left \lfloor \frac{n+1}{2}\right \rfloor $, then
	\[ \rho(G^{\prime})-\overline{ l}(G^{\prime}) > \rho(G)-\overline{ l}(G). \]
\end{lemma}

\begin{lemma}[\cite{wu2014}]\label{lemma 4 wu2014}
	Let $ G $ be a connected graph of order $ n\geq 3. $ Let $ u $ be a peripheral vertex of $ G $ and $ e_{G}(v) = d. $ Assume that $ G[V_{j}\cup V_{j+1}] $ is a clique for each $ j\in \{0,1,\dots,d-1\} $ and that $ n_{i}(u)\geq 2 $ for some $ i\in \{1,2,\dots,d-1\} $. Further, assume that $ i $ is the minimum subject to the above condition. Let $ v $ be a vertex in $ V_{i}(u) $ and $ V_{i-1}=\{u_{i-1}\} $. Let $ G^{\prime} $ be the graph with $ V (G^{\prime}) = V (G) $ and $ E(G^{\prime})=E(G)\cup A\setminus \{vu_{i-1}\}, $ where $ A = \{vy: y\in V_{i+2}\} $ if $ i\leq d-2 $, and $ A=\emptyset $ otherwise. If $ d = \left \lfloor \frac{n+1}{2}\right \rfloor, $ then
	\[ \rho(G^{\prime})-\overline{ l}(G^{\prime}) > \rho(G)-\overline{ l}(G). \]
\end{lemma}
A {\it Solt\'es} or {\it path-complete} graph \cite{soltes} is the graph obtained from a clique and a path by adding at least one edge between an endpoint of the path and the clique. The Solt\'es graphs are known to maximize the average distance $\overline{ l}$ when the number of vertices and of edges are fixed \cite{soltes}. 

Lemmas \ref{lemma 1 wu2014}, \ref{lemma 2 wu2014}, \ref{lemma 3 wu2014} and \ref{lemma 4 wu2014} lead the following result.

\begin{theorem}[\cite{wu2014}]\label{Thm 1 wu2014}
	Among all connected graphs $G$ on $n \ge 3$ vertices with average distance $\overline{ l}$ and remoteness $\rho$, the Solt\'es graphs with diameter $\lfloor (n+1)/2 \rfloor$
	maximize the difference $\rho - \overline{ l}$.
\end{theorem}

Theorem \ref{Thm 1 wu2014} can be equivalently stated as:
\begin{theorem}[\cite{wu2014}]
	Among all connected graphs $ G $ on $ n\geq 3 $ vertices with average distance $ \overline{ l} $ and remoteness $ \rho $, the maximum value of $ \rho-\overline{ l} $ is attained by the Solt\'es graphs with diameter $ D $, where 
	\[ \begin{cases}
		D=\frac{n+1}{2}, & \text{if $ n $ is odd}\\
		D\in \{\frac{n}{2}, \frac{n}{2}+1\} & \text{if $ n $ is even}.
	\end{cases} \]
\end{theorem}

Wu and Zhang \cite{wu2014} proved the following lemma, which helped them in proving Theorem \ref{Thm 2 wu2014}, earlier conjectured in \cite{Aouchiche2011}.

\begin{corollary}[\cite{wu2014}]
	Let $ G $ be connected graph with order $ n\geq 5 $. If $ n $ is odd and  $r=\frac{n-1}{2} $. Then $ \rho\geq \dfrac{n+1}{4} $, with equality if and only if $ G $ is the cycle $C_n$ or the graph composed of the cycle $C_n$ together with two crossed edges on four successive vertices of the cycle.
\end{corollary}
\begin{theorem}[\cite{wu2014}]\label{Thm 2 wu2014}
	Let $G$ be a connected graph on $n \ge 3$ vertices with remoteness $\rho$ and radius $r$. Then
	$$
	\rho - r \ge 
	\left\{
	\begin{array}{lll}
		\frac{3-n}{4} & & \mbox{ if $n$ is odd } \\
		\frac{n^2}{4n-4} - \frac{n}{2} & & \mbox{ if $n$ is even.}
	\end{array}
	\right.
	$$
	The inequality is best possible as shown by the cycle $C_n$ if $n$ is even and by the graph composed of the cycle $C_n$ together with two crossed edges on four successive vertices of the cycle.
\end{theorem}

\begin{figure}[h]
	\centerline{\scalebox{.75}{\includegraphics{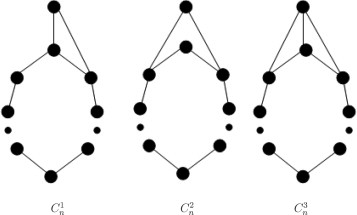}}}
	\caption{Graphs $ C_{n}^{1}, C_{n}^{2} $ and $ C_{n}^{1} $ (n is odd)}
	\label{Cn123}
\end{figure}
Hua, Chen and Das \cite{Hua2015} obtained the following result, earlier conjectured in \cite{Aouchiche2011}.
\begin{theorem}[\cite{Hua2015}]
	Let $ G $ be a connected graph on $ n \geq3 $ vertices with remoteness $ \rho $ and radius $ r $. If $ n $ is odd, then $ \rho - r \geq \frac{3-n}{4} $ with equality if and only if $ G\cong Cn $ or $ C_{n}^{i},~ i =1, 2, 3 $ (see Fig. \ref{Cn123}), and if $ n $ is even, then $ \rho-r \geq \frac{2n-n^{2}}{4(n-1)} $ with equality if and only if $ G\cong Cn. $
\end{theorem}

Next, we have result given upper bound for $ \rho-\pi $ in terms of order $ n $ and minimum degree $ \delta. $
\begin{theorem}[\cite{dankelmann2015}]\
	Let $ G $ be a connected graph of order $ n $ and minimum degree $ \delta $, where $ \delta\geq 2 $. Then
	\[\rho- \pi \leq \frac{3}{4(\delta + 1)}+ 3. \]
\end{theorem}

Dankelmann \cite{dankelmann2016} obtained some new bounds on proximity and remoteness.

\begin{theorem}[\cite{dankelmann2016}]
	Let $ G $ be a connected graph of order $ n $ and minimum degree $ \delta $, where $n\geq 20 $ and $ \delta\geq 2 $. Then
	\[ D-\pi \leq \frac{9n}{4(\delta+1)}+\frac{3\delta}{4}. \]
\end{theorem}
Next result gives sharp bound for remoteness in terms of diameter and order od graph
\begin{proposition}[\cite{dankelmann2016}]
	Let $ G $ be a connected graph of order $ n $ and diameter $ D $. Then
	\[ \rho \geq \frac{nd}{2(n-1)}, \]
	and this bound is sharp for all $ n $ and $ D $ with $ n\geq D + 1\geq 33 $ for which $ nD $ is even.
\end{proposition}

\begin{corollary}[\cite{dankelmann2016}]
	Let $ G $ be a connected graph of diameter $ D. $ Then 
	\[ \rho >\frac{D}{2}, \]
	and the coefficient $ \frac{1}{2} $ is best possible.
\end{corollary}

The next results gives the upper bound for $r-\pi$ in terms of minimum degree and order $n$ of $G.$
\begin{theorem}[\cite{dankelmann2016}]
	Let $ G$ be a connected graph of order $n$ and minimum degree $\delta$, where $ \delta<\frac{n}{4}-1. $ Then
	\[ r-\pi\leq \frac{3n}{4(\delta+1)} +\frac{8\delta+5}{4(\delta+1)}, \]
	and this bound is best possible, apart from an additive constant.
\end{theorem}
An immediate consequence of the above result is:
\begin{corollary}[\cite{dankelmann2016}]
	Let $ G $ be a connected graph. Then
	\[ \rho>\frac{r}{2}. \]
\end{corollary}

\begin{figure}[H]
	\begin{minipage}{.01\textwidth}
		~
	\end{minipage}
	\begin{minipage}{.45\textwidth}
		\centerline{\scalebox{.50}{\includegraphics{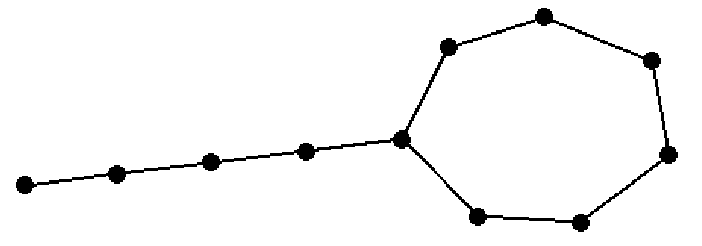}}}
		\caption{The lollipop $L_{11,7}$.}
		\label{lollipop}
	\end{minipage}
	\begin{minipage}{.02\textwidth}
		~
	\end{minipage}
	\begin{minipage}{.45\textwidth}
		\centerline{\scalebox{.50}{\includegraphics{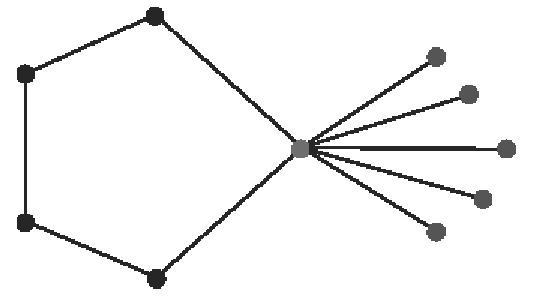}}}
		\caption{The turnip $T_{10,5}$.}
		\label{turnip}
	\end{minipage}
\end{figure}

A {\it lollipop } $L_{n,g}$ is the graph obtained from a cycle $C_g$ and a path $P_{n-g}$ by adding an edge between an endpoint of $P_{n-g}$ and a vertex of the cycle $C_g$. The lollipop $L_{11,7}$ is illustrated in Figure~\ref{lollipop}. For a lollipop $L_{n,g}$, we have
$$
\rho(L_{n,g}) = 
\left\{
\begin{array}{ccc}
	\frac{n}{2} - \frac{g(g-2)}{4(n-1)} & & \mbox{ if $g$ is even } \\
	\frac{n}{2} - \frac{g(g-2)+1}{4(n-1)} & & \mbox{ if $g$ is odd. }
\end{array}
\right.
$$
A {\it turnip} $T_{n,g}$, with $n \ge g \ge 3$, is the graph obtained from a cycle $C_g$ by attaching $n-p$ pending edges to one vertex from the cycle. The turnip $T_{10,5}$ is illustrated in Figure~\ref{turnip}. If $g = n$, the turnip $T_{n,g} = T_{n,n}$ is the cycle $C_n$. For a turnip $T_{n,g}$, we have
$$
\pi(T_{n,g}) = \left\{ 
\begin{array}{ccl}
	\frac{g^2-4g+4n-1}{4(n-1)} & & \mbox{ if $g$ is odd} \\
	\frac{g^2-4g+4n}{4(n-1)} & & \mbox{ if $g$ is even.}
\end{array}
\right.
$$

We have results from Aouchiche and Hensen\cite{Aouchiche2017} which relates $\pi$ with the girth of graphs. 

\begin{lemma}[\cite{Aouchiche2017}]\label{lemgirth}
	Let $G$ be a connected graph on $n \ge 3$ vertices with girth $g$. \\
	$(i)$ If $g = 3$, then $\pi \ge 1$ with equality if and only if $G$ contains a dominating vertex. \\
	$(ii) $ If $g = 4$, then $\pi \ge n/(n-1)$ with equality if and only if $G$ contains the turnip $T_{n,4}$ as a spanning subgraph and is a spanning subgraph of the complete bipartite graph $K_{n-2,2}$. \\
	$(iii) $ If $g \ge 5$, then $\pi \ge \pi(T_{n,g})$ with equality if and only if $G$ is the turnip $T_{n,g}$.
\end{lemma}

\begin{theorem}[\cite{Aouchiche2017}]\label{proxth}
	For any connected graph $G $ on $n \ge 3$ vertices with a finite girth $g $ and
	proximity $\pi $, we have
	\begin{equation}\label{pi1}
		\left.
		\begin{array}{lr} 
			\mbox{if $n $ is odd, } & \frac{1 - 3 n}{4}   \\
			\mbox{if $n $ is even, } & \frac{4 n - 3 n^{2}}{4 n - 4}
		\end{array}
		\right\} \le \pi - g \le \left\{
		\begin{array}{ll}
			\frac{n - 11}{4} - \frac{1}{n - 1} & \mbox{if $n $ is odd, } \\
			\frac{n - 11}{4} - \frac{3}{4(n - 1)} & \mbox{if $n $ is even; } 
		\end{array}
		\right.
	\end{equation}
	\begin{equation}\label{pi2}
		4 \le \pi + g \le \left\{
		\begin{array}{lr} 
			\frac{5n+1}{4} & \mbox{if $n $ is odd, }    \\
			\frac{5n^2 - 4n}{4(n-1)} & \mbox{if $n $ is even; } 
		\end{array}
		\right.
	\end{equation}
	\begin{equation}\label{pi3}
		\frac{1}{2\left\lfloor \sqrt{n} \right\rfloor + 1} + \frac{\left\lfloor \sqrt{n} \right\rfloor (\left\lfloor \sqrt{n} \right\rfloor -1)}{(2\left\lfloor \sqrt{n} \right\rfloor + 1)(n-1)} \le \frac{\pi}{g} \le 
		\left\{
		\begin{array}{ll} 
			\frac{n^{2} - 4}{12 n - 12} & \mbox{if $n $ is even, } \\
			\frac{n +1}{12} - \frac{1}{3 n - 3} & \mbox{if $n $ is odd; } 
		\end{array}
		\right.
	\end{equation}
	\begin{equation}\label{pi4}
		3 \le \pi \cdot g \le \left\{
		\begin{array}{lr} 
			\frac{n^2+n}{4} & \mbox{if $n $ is odd, }    \\
			\frac{n^3}{4(n-1)} & \mbox{if $n $ is even; } 
		\end{array}
		\right.
	\end{equation}
	The lower bound in (\ref{pi1}) and the upper bounds in (\ref{pi2}) and (\ref{pi4}) are reached if and only if $G$ is the cycle $C_n$. The upper bounds in (\ref{pi1}) and (\ref{pi3}) are reached if and only if $G$ is the lollipop $L_{n,3}$. The lower bounds in (\ref{pi2}) and (\ref{pi4}) are reached if and only if $G$ contains a dominating vertex. The lower bound in (\ref{pi3}) is reached if and only if $G$ is the turnip $T_{n, s}$, where $s = 2\left\lfloor \sqrt{n}\right\rfloor + 1$ when $\sqrt{n}$ is not an integer, and if and only if $G$ is any one of the turnips $T_{n,2\sqrt{n}-1 }$, $T_{n, 2\sqrt{n}}$ or $T_{n, 2\sqrt{n}+1}$ when $\sqrt{n}$ is an integer.
\end{theorem}

\begin{lemma}[\cite{Aouchiche2017}]\label{maxrho}
	Let $G$ be a connected graph on $n \ge 4$ vertices with girth $g \le n-1$ and remoteness $\rho$. Then 
	$$
	\rho \le \rho(L_{n,g})
	$$
	with equality if and only if $G$ is the lollipop $L_{n,g}$.
\end{lemma}

\begin{theorem}[\cite{Aouchiche2017}]
	For any connected graph $G $ on $n\ge 3$ vertices with remoteness $\rho $ and girth $g $, we have
	\begin{equation}\label{rho1}
		\left.
		\begin{array}{lr} 
			\mbox{if $n $ is even, } & \frac{4 n - 3 n^{2}}{4 n - 4}  \\
			\mbox{if $n $ is odd, }  & \frac{1 - 3 n}{4} 
		\end{array}
		\right\} \le \rho - g \le \frac{(n + 1)(n - 2)}{2 n - 2} -3 ;
	\end{equation}
	\begin{equation}\label{rho2}
		4 \le \rho + g \le \left\{
		\begin{array}{ll} 
			\frac{5 n^{2} - 4 n}{4 n - 4} & \mbox{if $n $ is even, } \\
			\frac{5 n +1}{4} & \mbox{if $n $ is odd; } 
		\end{array}
		\right.
	\end{equation}
	\begin{equation}\label{rho3}
		\frac{\rho}{g} \le \frac{(n + 1)(n - 2)}{6 n - 6};
	\end{equation}
	\begin{equation}\label{rho4}
		3 \le \rho \cdot g \le \rho(L_{n,g^*}) \cdot g^*
	\end{equation}
	where $g^*$ is the girth for which $\rho(L_{n,g_i}) \cdot g_i$, $i = 1 , \ldots, 4$, is maximum with $$g_1 = \left\lfloor \frac{2+\sqrt{6n^2-6n+4}}{3} \right\rfloor, g_2 = \left\lceil \frac{2+\sqrt{6n^2-6n+4}}{3} \right\rceil, g_3 = \left\lfloor \frac{2+\sqrt{6n^2-6n+7}}{3} \right\rfloor$$ and $g_4 = \left\lceil \frac{2+\sqrt{6n^2-6n+7}}{3} \right\rceil$.\\ 
	The lower bound in (\ref{rho1}) and the upper bound in (\ref{rho2}) are reached if and only if $G$ is the cycle $C_n$. The upper bounds in (\ref{rho1}) and (\ref{rho3}) are reached if and only if $G$ is the lollipop $L_{n,3}$. The lower bounds in (\ref{rho2}) and (\ref{rho4}) are reached if and only if $G$ is the complete graph $K_n$. The upper bound in (\ref{rho4}) is reached if and only if $G$ is the lollipop $L_{n,g^*}$.
\end{theorem}

The lower bound on the ratio of  $\rho$ and $g$ was first conjectured using AGX \cite{these},  and was later proved by Hua and Das in \cite{Hua2014}. 

\begin{theorem}[\cite{Hua2014}] \label{thm 1 hau2014}
	For any connected graph $G $ on $n\ge 3$ vertices with remoteness $\rho $ and girth $g $,
	$$
	\frac{\rho}{g} \ge 
	\left\{
	\begin{array}{ll} 
		\frac{n}{4 n - 4} & \mbox{if $n $ is even, } \\
		\frac{n + 1}{4 n} & \mbox{if $n $ is odd. } 
	\end{array}
	\right.
	$$
	with equality if and only if $G $ is a cycle $C_n $.
\end{theorem}

\begin{theorem}[\cite{Hua2014}] \label{thm 2 hau2014}
	Let $ G $ be a connected graph on $ n \geq 2 $ vertices with proximity $ \pi $ and average distance $ \overline{ l} $. Then 
	\[ \frac{\pi}{\overline{ l}}\geq \frac{n}{2(n-1)}, \]
	with equality if and only if $ G $ is isomorphic to the star $ S_{n }$.
\end{theorem}
The following is an immediate consequence of Theorem \ref{thm 2 hau2014}.
\begin{corollary}[\cite{Hua2014}] 
	Let $ G $ be a connected graph on $ n\geq 2 $ vertices with proximity $ \pi $, average degree $ \overline{d} $ and average distance $ \overline{ l} $. Then
	\[ \pi. \overline{d} \geq \overline{ l},\]
	with equality if and only if $ G $ is isomorphic to the star $ S_{n} $.
\end{corollary}

Dankelmann and Mafunda \cite{dankelmann2021a} gave results relating $ \pi $ with diameter and radius of triangle-free and $ C_{4} $-free graphs.

\begin{theorem}[\cite{dankelmann2021a}]
	If $ G $ is a connected, triangle-free graph of order $ n\geq 8$, minimum degree $ \delta\geq 3 $ and diameter $ D $, then
	\[ \pi \geq \frac{\delta(D-4)(D-1)}{8(n-1)}. \]
\end{theorem}

\begin{corollary}[\cite{dankelmann2021a}]
	If $ G $ is a connected, triangle-free graph of order $ n \geq 8$ and minimum degree $ \delta\geq 3 $, then
	\[ D-\pi \leq \frac{3(n-1)}{2\delta}+\frac{5}{2}. \]
	This bound is sharp apart from an additive constant.
\end{corollary}

\begin{theorem}[\cite{dankelmann2021a}]
	If $ G $ is a connected, triangle-free graph of order $ n \geq 6$, minimum degree $ \delta\geq 3 $ and radius $ r\geq 1 $, then
	\[ \pi \geq \frac{\delta}{2(n-1)} \left (r^{2}-7r+\frac{47}{8}\right) . \]
\end{theorem}

\begin{corollary}[\cite{dankelmann2021a}]
	If $ G $ is a connected, triangle-free graph of order $ n \geq 6$ and minimum degree $ \delta\geq 3 $, then
	\[ r-\pi \leq \frac{n-1}{2\delta}+\frac{11}{2}. \]
	This bound is sharp apart from an additive constant.
\end{corollary}

\begin{theorem}[\cite{dankelmann2021a}]
	If $ G $ is a connected, $ C_{4} $-free graph of order $ n \geq 16$ and minimum degree $ \delta\geq 3 $, then
	\[ \pi \geq \frac{\delta^{2}-2\left \lfloor \frac{\delta}{2}\right \rfloor+1}{5(n-1)}\left(r^{2}-8r+\frac{127}{8}\right) . \]
\end{theorem}

\begin{corollary}[\cite{dankelmann2021a}]
	If $ G $ is a connected, $ C_{4} $-free graph of order $ n \geq 16$ and minimum degree $ \delta\geq 3 $, then
	\[r- \pi \leq \frac{5(n-1)}{4\left(\delta^{2}-2\left \lfloor \frac{\delta}{2}\right \rfloor+1\right) }+4. \]
\end{corollary}

Czabarka et al. \cite{czabarka2021} gave the upper bounds on $ \pi $ in triangulations and quadrangulations.

The plane graph $ G $ is a triangulation (respectively quadrangulation) if every face is a triangle (respectively 4-cycle).
\begin{proposition}[\cite{czabarka2021}]
	\begin{itemize}
		\item [(a)]
		Let $ G $ be a $ 5 $-connected triangulation of order $ n $. Then
		\[ \rho\leq \frac{n + 4}{10}+\epsilon_{ n}, \]
		where $\epsilon_{ n}= -\frac{3}{5(n-1)}$ if $ n \equiv 0 $ (mod $ 5 $), $\epsilon_{ n}= - \dfrac{1}{n-1} $ if $ n\equiv 1 $ (mod $ 5 $), $\epsilon_{ n} = \dfrac{2}{5(n-1)} $ if $ n\equiv 2 $ (mod $ 5 $), and $\epsilon_{ n} = - \frac{2}{5(n-1)} $ if $ n \equiv 3, 4 $ (mod $ 5 $).
		\item [(b)] If $ G $ is a $ 3 $-connected quadrangulation of order $ n $, then
		\[ \rho \leq \frac{n+2}{6}+\epsilon_{ n}, \]
		where $ \epsilon_{ n} = -\frac{5}{3(n-1)} $ if $ n\equiv 0$ (mod $ 3 $), $ \epsilon_{ n} = - \frac{1}{n-1}  $ if $ n \equiv 1 $ (mod $ 3 $), and $ \epsilon_{ n} = \dfrac{1}{3(n-1)} $ if $ n \equiv 2 $(mod $ 3 $).
	\end{itemize}
\end{proposition}

\section{Proximity and Remoteness Compared to other Invariants}
Aouchiche and Hasen \cite{Aouchiche2011} related remoteness with the independence number $\alpha$ and  matching number $\mu$ of $G$ and presented the following results.

\begin{proposition}[\cite{Aouchiche2011}]\label{prop-alpha}    
	Let $G$ be a connected graph on $n \ge 8$ vertices with remoteness $\rho $ and independence number
	$\alpha$. Then 
	$$
	\rho - \alpha \ge 3-n-\frac{1}{n-1} 
	$$
	and 
	$$
	\rho - \alpha \le \left\{
	\begin{array}{lcl}
		\frac{n-3}{8} - \frac{3}{8(n-1)}  &  & \mbox{ if  }  n = 0 \, \, (mod \, \, 4), \\
		\frac{n-3}{8}   &  & \mbox{ if  }  n = 1 \, \, (mod \, \, 2), \\
		\frac{n-3}{8} + \frac{1}{8(n-1)}  &  & \mbox{ if  }  n = 2 \, \, (mod \, \, 4). 
	\end{array}
	\right. 
	$$
	The lower bound is attained if and only if $G$ is the star $S_n$. The bound is best possible as shown by the kites $KI_{n,n_0}$ where $n_0 = (n+ n \,\, (mod \, \,  4))/2$.
\end{proposition}

\begin{proposition}[\cite{Aouchiche2011}]   
	Let $G$ be a connected graph on $n \ge 3$ with remoteness $\rho$ and matching number $\mu$. Then
	$$
	\rho - \mu \le \left\{
	\begin{array}{lll}
		\frac{n+1}{8} + \frac{1}{8(n-1)} & & \mbox{ if } n = 0 \, \, (mod \, \,  4), \\
		\frac{n+1}{8}  & & \mbox{ if } n = 1 \, \, (mod \, \,  4), \\
		\frac{n+1}{8} - \frac{3}{8(n-1)} & & \mbox{ if } n = 2 \, \, (mod \, \,  4), \\
		\frac{n+1}{8}  & & \mbox{ if } n = 3 \, \, (mod \, \,  4), \\
	\end{array}
	\right.
	$$
	with equality if and only if $G$ is the comet $CO_{n,n-D^*+1}$ where $D^* = 2\left\lfloor \frac{n+2}{4} \right\rfloor$.
\end{proposition}

The ({\it vertex}) {\it connectivity} $\nu$ of a connected graph $G$ is the minimum number of vertices whose removal disconnects $G$ or reduces it to a single vertex. The {\it algebraic connectivity} $a$ of a graph $G$ is the second smallest eigenvalue of its Laplacian matrix $L = Diag - A$, where $Diag$ is the diagonal square matrix indexed by the vertices of $G$ whose diagonal entries are the degrees in $G$, and $A$ is the adjacency matrix of $G$.

Sedlar et al. \cite{sedlar2008} related vertex connectivity $\nu$, algebraic connectivity $a$ and remoteness $\rho$ in the following result.
\begin{theorem}[\cite{sedlar2008}]
	Let $G$ be a connected graph on $n \ge 2$ vertices with vertex connectivity $\nu$, algebraic connectivity $a$ and remoteness $\rho$. Then
	$$
	\nu \cdot \rho \le n-1
	$$
	with equality if and only if $G$ is the complete graph $K_n$; and 
	$$
	a \cdot \rho \le n
	$$
	with equality if and only if $G$ is the complete graph $K_n$. Moreover, if $G$ is not complete, then
	$$
	a \cdot \rho \le n-1 - \frac{1}{n-1}
	$$
	with equality if and only if $G \cong K_n - M$, where $M$ is any non empty set of disjoint edges.
\end{theorem}

The following sequence of results of Hua and Das \cite{Hua2014} presented results for $\rho$ and $\pi$ in terms of clique number, average degree and average distance of $G.$

\begin{theorem}[\cite{Hua2014}]
	Let $ G $ be a connected graph on $ n \geq 3 $ vertices with remoteness $ \rho $ and maximum degree $ \triangle $. If $ \triangle \geq \left \lceil \frac{n}{4} \right \rceil +1, $ then 
	\[ \rho+\triangle\geq \begin{cases}
		\frac{n+9}{4} & \text{if $ n $ is odd},\\
		2+\frac{n^{2}}{4n-4} & \text{if $ n $ is even},
	\end{cases} \]
	with equality if and only if $ G\cong C_{3} $ or $ G\cong C_{4}. $
\end{theorem}

\begin{theorem}[\cite{Hua2014}]
	Let $ G $ be a connected graph on $ n $ vertices, $ m $ edges with proximity $ \pi $ and average degree $ \overline{d} $. If $ m\leq \frac{2(n-1)^{2}}{n}, $ then 
	\[ \pi .\overline{d}\leq n-1. \]
\end{theorem}

\begin{theorem}[\cite{Hua2014}]
	Let $ G $ be a connected graph on $ n \geq 3 $ vertices with clique number $ \omega $ and remoteness $ \rho $. Then 
	\[ \rho \leq \frac{n^{2}-\omega^{2}-n+3\omega-2}{2(n-1)}, \]
	with equality if and only if $ G\cong Ki_{n,\omega}. $
\end{theorem}

\begin{theorem}[\cite{Hua2014}]
	Let $ G $ be a connected graph on $ n \geq 3 $ vertices with clique number $ \omega $, remoteness $ \rho $ and proximity $ \pi $. Then 
	\[ \rho+\pi \leq \begin{cases}
		\frac{3n^{2}-2\omega^{2}-2n+6\omega-5}{4(n-1)} & \text{if $ n $ is odd},\\
		\frac{3n^{2}-2\omega^{2}-2n+6\omega-4}{4(n-1)} & \text{if $ n $ is odd},
	\end{cases} \]
	with equality if and only if $ G\cong P_{n}. $
\end{theorem}
The next result gives an upper bound on $ \rho-\pi $ in terms of $ n $ and $ \omega. $
\begin{theorem}[\cite{Hua2014}]
	Let $ G $ be a connected graph on $ n \geq 3 $ vertices with clique number $ \omega $, remoteness $ \rho $ and proximity $ \pi $. Then 
	\[ \rho-\pi \leq \frac{(n-\omega)(n+\omega-3)}{2(n-1)},
	\]
	with equality if and only if $ G\cong Ki_{n,n-1}. $
\end{theorem}

\begin{theorem}[\cite{Hua2014}]
	Let $ G $ be a connected bipartite graph with each partite set of cardinality $ n\geq 2 $. If the average distance 
	\[ \overline{ l}\leq \frac{3}{2}+\frac{n-2}{2n(2n-1)}, \]
	then $ G $ is Hamiltonion.
\end{theorem}

\begin{corollary}[\cite{Hua2014}]
	Let $ G $ be a connected bipartite graph with each partite set of cardinality $ n\geq 2 $. If 
	\[ \rho\leq\frac{3}{2}+\frac{n-2}{2n(2n-1)}, \]
	then $ G $ is Hamiltonion.
\end{corollary}

The distance eigenvalues of a connected graph, denoted by $\partial_1, \partial_2, \ldots, \partial_n$, are those of its distance matrix, and are indexed such that $\partial_1 \ge \partial_2 \ge \ldots \ge \partial_n$. For a detailed survey on distance spectra of graphs see \cite{distsurvey}. Next, we have results relating proximity $\pi$ and remoteness $\rho$ with the distance eigenvalues of a connected graph $G$.

\begin{theorem}[\cite{Aouchiche2016}]\label{minmax}
	Let $G$ be a graph on $n \ge 4$ vertices with largest distance eigenvalue $\partial_1$, proximity $\pi$ and remoteness $\rho$. Then 
	$$
	\pi \le \overline{ l} \le \frac{\partial_1}{ n-1} \le \rho
	$$
	with equalities if and only if $G$ is a transmission regular graph.
\end{theorem}

\begin{corollary}[\cite{Aouchiche2016}]
	Let $G$ be a graph on $n \ge 2$ vertices with largest distance eigenvalue $\partial_1$ and proximity $\pi$. Then 
	$$
	\partial_1 - \pi \ge n-2
	$$
	with equalities if and only if $G$ is the complete graph $K_n$.
\end{corollary}

\begin{corollary}[\cite{Aouchiche2016}]
	Let $G$ be a graph on $n \ge 4$ vertices with second largest distance eigenvalue $\partial_2$ and remoteness $\rho$. Then 
	$$
	\rho + \partial_2 \ge 0
	$$
	with equality if and only if $G$ is the complete graph $K_n$.
\end{corollary}
The bound in the above corollary is best possible among the bounds of the form $\rho + \partial_k \ge 0$, with a fixed integer $k$, over the class of all connected graphs. Indeed, if we consider the complete bipartite graphs $K_{\left\lfloor n/2 \right\rfloor, \left\lceil n/2 \right\rceil}$, on $n \ge 3$, by direct calculation, we get
$$
\rho(K_{\left\lfloor n/2 \right\rfloor, \left\lceil n/2 \right\rceil}) + \partial_3(K_{\left\lfloor n/2 \right\rfloor, \left\lceil n/2 \right\rceil}) = \left\{ 
\begin{array}{lcl}
	-\frac{1}{2} &  & \mbox{ if $n$ is odd,} \\
	\\
	-\frac{1}{2} - \frac{1}{2(n-1)}  &  & \mbox{ if $n$ is even,} \\
\end{array}
\right.
$$
which is negative for $n \ge 3$.

\begin{proposition}[\cite{Aouchiche2016}]
	Let $T$ be a tree on $n \ge 4$ vertices with remoteness $\rho$, diameter $D$ and distance spectrum $\partial_1 \ge \partial_2 \ge \cdots \ge \partial_n$. Then 
	$$
	\rho + \partial_{\left\lfloor \frac{D}{2} \right\rfloor} > 0.
	$$
\end{proposition}

\begin{corollary}[\cite{Aouchiche2016}]\label{cor-3}
	Let $G$ be a graph on $n \ge 4$ vertices with second largest distance eigenvalue $\partial_2$ and proximity $\pi$. Then 
	$$
	\pi + \partial_n \le 0
	$$
	with equality if and only if $G$ is the complete graph $K_n$.
\end{corollary}

\begin{theorem}[\cite{Aouchiche2016}]
	Let $G$ be a graph on $n \ge 4$ vertices with largest distance eigenvalue $\partial_1$ and remoteness $\rho$. Then 
	$$
	\partial_1 - \rho \ge n-2
	$$
	with equalities if and only if $G$ is the complete graph $K_n$.
\end{theorem}

\begin{proposition}[\cite{Aouchiche2016}]\label{prop-rho}
	Let $G$ be a graph on $n \ge 4$ vertices with least distance eigenvalue $\partial_n$ and remoteness $\rho$. Then   
	$$
	\partial_n + \rho \le 0 
	$$
	with equality if and only if $G$ is $K_n$.
\end{proposition}

\begin{theorem}[\cite{Aouchiche2016}]\label{cor-pi+indice}
	Let $G$ be a graph on $n \ge 4$ vertices with second largest distance eigenvalue $\partial_2$ and proximity $\pi$. Then 
	$$
	\pi + \partial_2 \ge 0
	$$
	with equality if and only if $G$ is the complete graph $K_n$.
\end{theorem}

\begin{proposition}[\cite{Aouchiche2016}]
	Let $T$ be a graph on $n \ge 4$ vertices with proximity $\pi$, diameter $D$ and distance spectrum $\partial_1 \ge \partial_2 \ge \cdots \ge \partial_n$. Then 
	$$
	\pi + \partial_{\left\lfloor \frac{D}{2} \right\rfloor} > 0.
	$$
\end{proposition}

Lin, Das and Wu \cite{Lin2016} proved two theorems, earlier as conjectures in \cite{Aouchiche2016}. Following lemmas and results were proved in \cite{Lin2016}.
\begin{lemma}[\cite{Lin2016}]
	Let $ G $ be a connected graph of order $ n $ with diameter $ D $ and remoteness $ \rho $. Then
	\[ \rho \geq \frac{D}{2}. \]
\end{lemma}

\begin{theorem}[\cite{Lin2016}]
	Let $ G $ be a connected graph of order $ n \geq 4 $ with diameter $ D $, remoteness $ \rho $ and distance eigenvalues $ \partial_{1} \geq \dots \geq  \partial_{n}. $ Then we have the following statements.
	\begin{itemize}
		\item[\bf (i)] If $ D  = 2,$ then
		\[ \rho + \partial_3 \geq \frac{\left \lceil \frac{n}{2}\right \rceil-2}{n-1}-1, \]with equality holding if and only if $ G\cong K_{n_{1},n_{2}}. $
		\item[\bf (ii)] If $ D \geq  3 $, then
		\[ \rho + \partial_3 >\frac{d}{2}- 1.2. \]
	\end{itemize}
\end{theorem}

\begin{theorem}[\cite{Lin2016}]
	Let $ G $ be a connected graph of order $ n \geq  4 $ with diameter $ D $, remoteness $ \rho $ and distance spectrum $ \partial_{1} \geq \dots \geq  \partial_{n}. $ Then
	\[  \rho + \partial_{ \left \lfloor \frac{7D}{8}\right   \rfloor} >0.  \]
\end{theorem}

Besides the above results, more results regarding remoteness and distance eigenvalues were given in the same article \cite{Lin2016}. Before stating them, we need the following definition.

Denote by $ H_{n-D} (n>D) $ a graph of order $ n - D $ such that $ V(H_{n-d}) = V(\overline{K}_{n-d}) $ and $ E(H_{n-D}) \supseteq E(\overline{K}_{n-D}), $ where $\overline{K}_{n-D}$ is a null graph of order $ n-D. $ Let $ H_{n,D} $ be a graph of order $ n $ with diameter $ D $ obtained by joining $ n-D $ edges between one end of the path $ P_{D} $ with each vertex of $ H_{n-D}. $ Now, we are in a position to state the remaining results  from \cite{Lin2016}.

\begin{lemma}[\cite{Lin2016}]
	Let $ G $ be a connected graph of order $ n $ with diameter $ D $ and remoteness $ \rho $. Then
	\[ \rho \leq  D -\frac{D^{2} -D}{2(n - 1)}, \]
	with equality holding if and only if $ G\cong H_{n, D}. $
\end{lemma}

\begin{theorem}[\cite{Lin2016}]
	Let $ G $ be a connected graph of order $ n $ with diameter $ D $ and remoteness $ \rho $. Then
	\[ \rho +\partial_n \leq -\frac{D^{2} -D}{2(n - 1)}, \]
	with equality holding if and only if $ G\cong K_{n}. $
\end{theorem}

\begin{lemma}[\cite{Lin2016}]
	Let $ G $ be a connected graph of order $ n $ with diameter $ D\geq 3. $ Then
	\[ \partial_{1}>n-2+D. \]
\end{lemma}

\begin{theorem}[\cite{Lin2016}]
	Let $ G \ncong K_{n}$ be a connected graph of order $ n $  and remoteness $ \rho $. Then
	\[ \partial_1-\rho \geq \frac{n-1+\sqrt{(n-1)^{2}+8}}{2}-\frac{n}{n-1}, \]
	with equality holding if and only if $ G\cong K_{n}-e, $ where $ e $ is an edge of $ G. $
\end{theorem}

Jia and Song \cite{jia2018} obtained various results related to remoteness, distance eigenvalues, distance (signless) Laplacian eigenvalues of graphs.
\begin{theorem}[\cite{jia2018}]
	Let $ G $ be a connected graph of order $ n \geq 4 $ with remoteness $ \rho $. Then
	\[ n\leq \rho + \partial_{1} \leq  \rho(P_{n}) + \partial_{1}(P_{n}), \]
	with the left equality holding if and only if $ G\cong K_{n} $ and the right equality holding if and only if $ G\cong P_{n}. $
\end{theorem}

\begin{theorem}[\cite{jia2018}]
	Let $ G\ncong K_{n} $ be a connected graph of order $ n \geq 4 $ with remoteness $ \rho $. Then
	\[ \rho+\partial_{1}\geq \frac{n}{n-1}+\frac{n -1+\sqrt{(n -1)^{2} + 8}}{2},\]
	with equality holding if and only if $ G\cong K_{n}-e. $
\end{theorem}

\begin{theorem}[\cite{jia2018}]
	Let $ G\ncong (K_{n}, K_{n}-e) $ be a connected graph of order $ n \geq 4 $ with remoteness $ \rho $. Then
	\[ \rho+\partial_{1}\geq \frac{n}{n-1}+\frac{n -1+\sqrt{(n -1)^{2} + 16}}{2},\]
	with equality holding if and only if $ G\cong K_{n}-2e, $ where $ 2e $ are two matching edges.
\end{theorem}

\begin{theorem}[\cite{jia2018}]
	Let $ G $ be a complete bipartite graph of order $ n \geq 4 $ with remoteness $ \rho $. Then
	\[ \rho+\partial_2\geq n-\frac{1}{n-1}-\sqrt{n^{2}-3n+3}, \]
	with equality holding if and only if $ G $ is star.
\end{theorem}

For connected graphs Jia and Song \cite{jia2018} proposed the following conjecture.
\begin{conjecture}[\cite{jia2018}]
	Let $ G\ncong (K_{n}, K_{n}-e) $ be a connected graph of order $ n \geq 4 $ with remoteness $ \rho $. Then
	\[ \rho+\partial_{2}\geq \frac{n}{n-1}+\frac{n -1-\sqrt{(n -1)^{2} + 8}}{2},\]
	with equality holding if and only if $ G\cong K_{n}-2e, $ where $ 2e $ are two matching edges.
\end{conjecture}

Other results of Jia and Song \cite{jia2018} are stated as:
\begin{theorem}[\cite{jia2018}]
	Let $ G\ncong K_{n} $ be a connected graph of order $ n\geq 4 $ with remoteness $ \rho $. Then
	\[\rho-\partial_{n} \geq 2, \]
	with equality holding if and only if $ G\cong K_{n}. $
\end{theorem}

\begin{theorem}[\cite{jia2018}]
	Let $ G\ncong K_{n} $ be a connected graph of order $ n\geq 4 $ with remoteness $ \rho $. Then
	\[\rho-\partial_{n} \geq 3+\frac{1}{n-1}, \]
	with equality holding if and only if $ G\cong K_{n}-me, $ where $ me $ denotes $ m $ matching edges.
\end{theorem}

The distance Laplacian eigenvalues of a connected graph, denoted by $\partial_{1}^{L}, \partial_{2}^{L}, \ldots, \partial_{n}^{L}$, are the eigenvalues of its distance Laplacian matrix, and are indexed such that $\partial_{1}^{L}\geq \partial_{2}^{L}\geq  \ldots \geq \partial_{n}^{L}$, where $\partial_{1}^{L}$ is known as the distance Laplacian spectral radius of $G.$ The following results of \cite{jia2018} presents bounds for $\rho,\partial_{1}^{L}$ and $\partial_{1}^{Q}$ and identifies the candidate graphs attaining the equalities. 

\begin{theorem}[\cite{jia2018}]
	Let $ G $ be a connected graph of order $ n\geq 4 $ with remoteness $ \rho $. Then
	\[ n+1\leq \rho +\partial_{1}^{L}\leq \rho(P_{n})+\partial_{1}^{L}(P_{n}), \]
	with the left equality holding if and only if $ G\cong K_{n} $ and the right equality holding if and only if $ G\cong P_{n}. $
\end{theorem}

\begin{theorem}[\cite{jia2018}]
	Let $ G\ncong K_{n} $ be a connected graph of order $ n\geq 4 $ with remoteness $ \rho $. Then
	\[ \rho +\partial_{1}^{L}\geq n+\frac{1}{n-1}+3, \]
	with equality holding if and only if $ G\cong K_{n}-me, $ where $ me $ denotes $ m $ matchings.
\end{theorem}

\begin{theorem}[\cite{jia2018}]
	Let $ G $ be a connected graph of order $ n\geq 4 $ with remoteness $ \rho $. Then
	\[ \partial_{1}^{L}-\rho \geq n-1, \]
	with equality holding if and only if $ G\cong K_{n}. $
\end{theorem}

\begin{theorem}[\cite{jia2018}]
	Let $ G \ncong K_{n} $ be a connected graph of order $ n\geq 5 $ with remoteness $ \rho $. Then
	\[ \partial_{1}^{L}-\rho \geq n+1-\frac{1}{n-1}, \]
	with equality holding if and only if $ G\cong K_{n}-me, $ where $ 1\leq m\leq \left \lfloor \frac{n}{2}\right \rfloor. $
\end{theorem}

The distance signless Laplacian eigenvalues of a connected graph, denoted by $\partial_{1}^{Q}, \partial_{2}^{Q}, \ldots, \partial_{n}^{Q}$, are those of its distance signless Laplacian matrix, and are indexed such that $\partial_{1}^{Q}\geq \partial_{2}^{Q}\geq  \ldots \geq \partial_{n}^{Q}$.
\begin{theorem}[\cite{jia2018}]
	Let $ G $ be a connected graph of order $ n\geq 4 $ with remoteness $ \rho $. Then
	\[ 2n\leq 2\rho +\partial_{1}^{Q}\leq 2\rho(P_{n})+\partial_{1}^{Q}(P_{n}), \]
	with the left equality holding if and only if $ G\cong K_{n} $ and the right equality holding if and only if $ G\cong P_{n}. $
\end{theorem}

\begin{theorem}[\cite{jia2018}]
	Let $ G\ncong K_{n} $ be a connected graph of order $ n\geq 4 $ with remoteness $ \rho $. Then
	\[ 2\rho +\partial_{1}^{Q}\geq \frac{3n-2+\sqrt{(n-2)^{2}+16}}{2}+\frac{2n}{n-1}, \]
	with equality holding if and only if $ G\cong K_{n}-e. $ 
\end{theorem}

\begin{theorem}[\cite{jia2018}]
	Let $ G\ncong (K_{n}, K_{n}-e) $ be a connected graph of order $ n\geq 4 $ with remoteness $ \rho $. Then
	\[ 2\rho +\partial_{1}^{Q}\geq \frac{3n-2+\sqrt{(n-2)^{2}+32}}{2}+\frac{2n}{n-1}, \]
	with equality holding if and only if $ G\cong K_{n}-2e, $ where $ 2e $ are two matching edges. 
\end{theorem}

\begin{theorem}[\cite{jia2018}]
	Let $ G $ be a connected graph of order $ n\geq 4 $ with remoteness $ \rho. $ Then
	\[ \partial_{1}^{Q}-2\rho \geq 2n-4, \]
	with equality holding if and only if $ G\cong K_{n}. $
\end{theorem}

\begin{theorem}[\cite{jia2018}]
	Let $ G\ncong K_{n} $ be a connected graph of order $ n\geq 4 $ with remoteness $ \rho $. Then
	\[ \partial_{1}^{Q}-2\rho \geq \frac{3n-2+\sqrt{(n-2)^{2}+16}}{2}-\frac{2n}{n-1}, \]
	with equality holding if and only if $ G\cong K_{n}-e $. 
\end{theorem}

Mojallal and Hansen \cite{mojallal2021} obtained relation between proximity and the third largest distance eigenvalue $\partial_3(G)$ of a graph $G$, some results earlier conjectured in \cite{Aouchiche2016}.

\begin{figure}[h]
	\centerline{\scalebox{.75}{\includegraphics{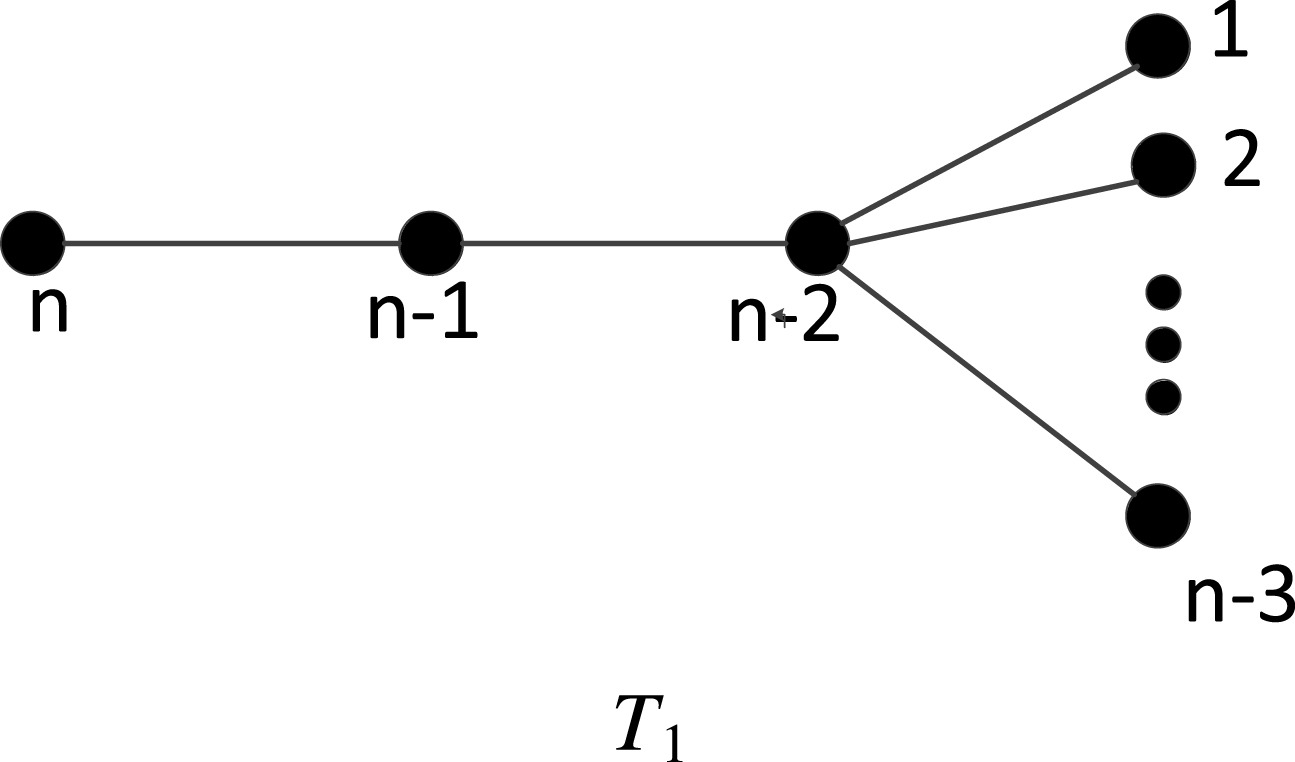}}\qquad \scalebox{.75}{\includegraphics{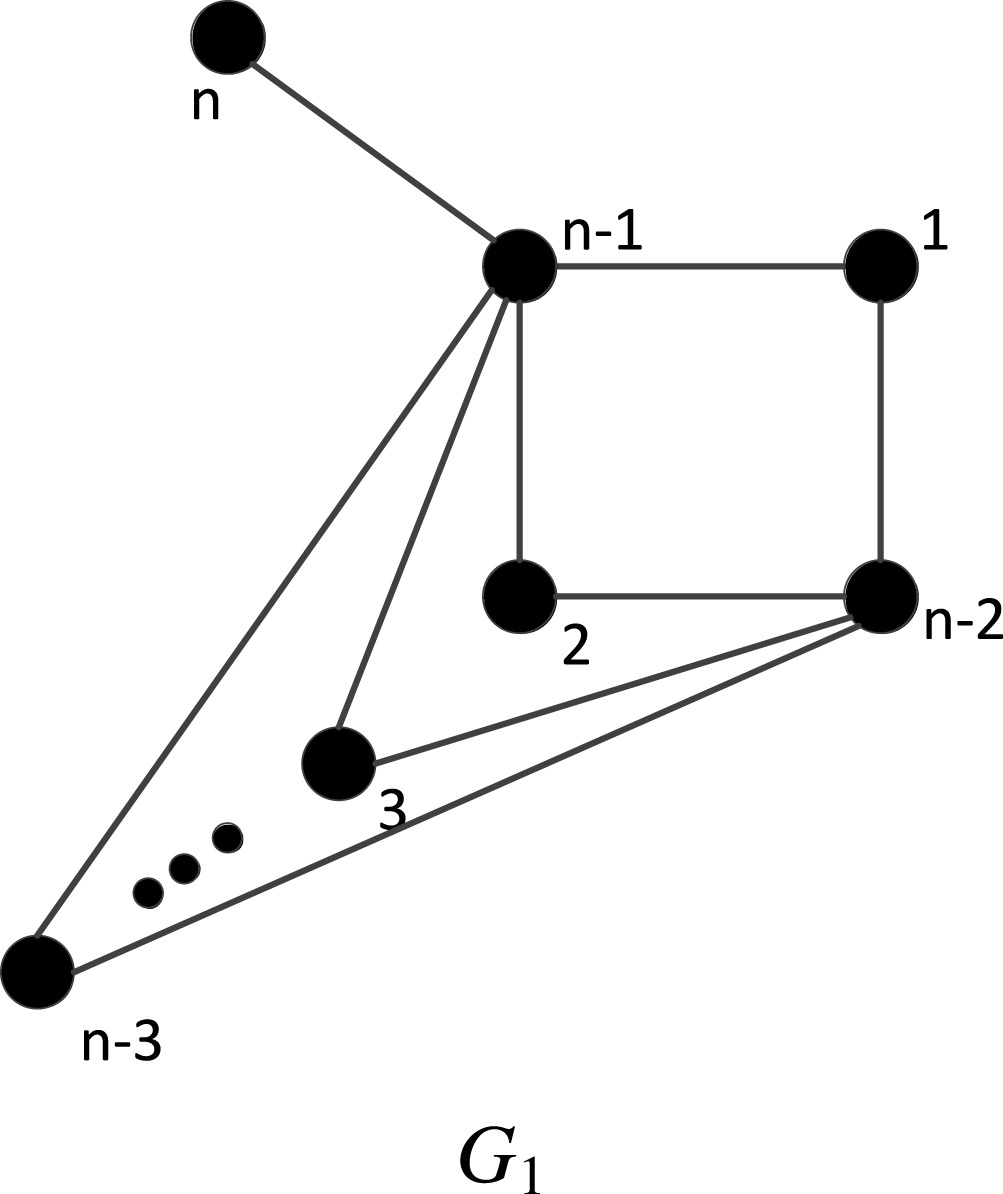}}}
	\caption{The tree $ T_{1} $ and the graph $ G_{1} $}
	\label{tree T1}
\end{figure}
\begin{lemma}[\cite{mojallal2021}]
	Let $ T_{1} $ be the tree of order $ n $ given in Fig \ref{tree T1}. Then 
	\[ \pi(T_{1})+ \partial_3(T_{1})>0. \]
\end{lemma}

\begin{lemma}[\cite{mojallal2021}]
	Let $ G_{1} $ be the graph of order $ n $ given in Fig \ref{tree T1}. Then
	\[ \pi((G_{1}))+\partial_3(G_{1})>0. \]
\end{lemma}

\begin{lemma}[\cite{mojallal2021}]
	Let $ G $ be a graph of order $ n $ with the diameter $ D = 3 $ and let $ i_{1} $ and $ i_{4} $ be two vertices of $ G $ with distance $ 3 $. If $ d_{G}(i_{1}) \geq 2 $ and $ d_{G}(i_{4}) \geq 2. $ Then 
	\[ \pi+ \partial_3 >0. \]
\end{lemma}

\begin{lemma}[\cite{mojallal2021}]
	Let $ G $ be a graph of order $ n $ with the diameter $ D=3 $ and let $ i_{1} $ and $ i_{4} $ be two vertices of $ G $ with distance $ 3 $. If $ d_{G}(i_{1}) = 1 $ or $ d_{G}(i_{4}) = 1. $ Then \[ \pi + \partial_3> 0. \]
\end{lemma}

\begin{theorem}[\cite{mojallal2021}]
	Let $ G $ be a graph with the diameter $ D\geq 3 $, proximity $ \pi(G) $ and third largest distance eigenvalue $ \partial_3. $ Then
	\[ \pi + \partial_3 > 0. \]
\end{theorem}

\begin{corollary}[\cite{mojallal2021}]
	Let $ G $ be a graph with the diameter $ D \geq 3, $ and the third largest distance eigenvalue $ \partial_3. $ Then
	\begin{align*}
		(i)~ D + \partial_3 > 0\quad (ii)~ \rho + \partial_3 > 0\quad (iii)~ \overline{ l}+ \partial_3 > 0 \quad (iv)~ ecc+ \partial_3 > 0\quad (v)~r + \partial_3> 0.
	\end{align*}
\end{corollary}

\section{Conclusion}
The survey collects all the results related to proximity $\pi$ and remoteness $\rho$ and its relation with other spectral invariants published in different journals from 2006 to January 2024.

\section*{Data Availability}
There is no data associated with this article.

\section*{Conflict of interest}
The authors declare that they have no competing interests.

\end{document}